\documentclass[11 pt]{amsart}

\usepackage[latin1]{inputenc}
\usepackage{amsmath,amsthm,amssymb,amscd,graphicx,psfrag}

\newtheorem{thm}{Theorem}[section]
 
 \newtheorem{lem}[thm]{Lemma}
 \newtheorem{prop}[thm]{Proposition}
 \theoremstyle{definition}
 \newtheorem{df}[thm]{Definition}
 \theoremstyle{remark}
 \newtheorem{rem}[thm]{Remark}
 \newtheorem{ex}{Example}
 \numberwithin{equation}{section}

\def\be#1 {\begin{equation} \label{#1}}
\newcommand{\ee}{\end{equation}}
\def\dem {\noindent \textsc{Proof:} }

\def\sqw{\hbox{\rlap{\leavevmode\raise.3ex\hbox{$\sqcap$}}$%
\sqcup$}}
\def\findem{\ifmmode\sqw\else{\ifhmode\unskip\fi\nobreak\hfil
\penalty50\hskip1em\null\nobreak\hfil\sqw
\parfillskip=0pt\finalhyphendemerits=0\endgraf}\fi}

\newcommand{\mb}{\medskip\noindent}
\newcommand{\gb}{\bigskip\noindent}
\newcommand{\R}{\mathbb R}

\newcommand{\N}{\mathbb N}
\newcommand{\Z}{\mathbb Z}

\newcommand{\s}{\mathcal S}

\newcommand\<{\langle}
\renewcommand\>{\rangle}

\begin{document}

\author{Fr\'ed\'eric Bernicot}
\address{CNRS - Universit\'e Lille 1 \\ Laboratoire de math\'ematiques Paul Painlev\'e \\ 59655 Villeneuve d'Ascq Cedex, France}
\curraddr{}
\email{frederic.bernicot@math.univ-lille1.fr}

\author{Pierre Germain}
\address{Courant Institute of Mathematical Sciences \\ New York University \\ 251 Mercer Street \\ New York, N.Y. 10012-1185 \\ USA}
\curraddr{ETH Z\"urich \\ Departement Mathematik \\ R\"amistrasse 101 \\ 8092 Zurich \\ Switzerland}
\email{pgermain@cims.nyu.edu}

\title{Bilinear oscillatory integrals and boundedness for new bilinear multipliers}

\subjclass[2000]{Primary 42B10 ; 42B20}

\keywords{Bilinear multipliers ; oscillatory integrals}

\begin{abstract} We consider bilinear oscillatory integrals, i.e. pseudo-product operators whose symbol involves an oscillating factor. Lebesgue space inequalities are established, which give decay as the oscillation becomes stronger ; this extends the well-known linear theory of oscillatory integral in some directions. The proof relies on a combination of time-frequency analysis of Coifman-Meyer type with stationary and non-stationary phase estimates. As a consequence of this analysis, we obtain Lebesgue estimates for new bilinear multipliers defined by non-smooth symbols.
\end{abstract}

\date{\today}

\maketitle
\tableofcontents

\section{Introduction}

\subsection{Presentation of the problem}

\subsubsection{Pseudo-products}

Pseudo-product operators were first introduced by R. Coifman and Y. Meyer~\cite{coifmanmeyer}. These are the multilinear operators mapping functions on $\mathbb{R}^d$ to a function on $\mathbb{R}^d$ which are invariant by space translation.
Turning for simplicity in the notations to the bilinear case, a pseudo-product operator can be written as
\begin{align}
T_{m}(f,g)(x) & := \mathcal{F}^{-1} \left[\int_{\mathbb{R}^d} m(\eta,\xi-\eta) \widehat{f}(\eta) \widehat{g} (\xi-\eta) \,d\eta \right](x) \label{defB} \\
 & = \int_{\mathbb{R}^{2d}} e^{ix\cdot(\eta+\xi)} m(\eta,\xi) \widehat{f}(\eta) \widehat{g} (\xi) \,d\eta d\xi  \,\,, \nonumber
\end{align}
where $m$ is the \textit{symbol} of the operator, and the Fourier transform of $f$ is denoted by $\mathcal{F} f$ or
$\widehat{f}$ (see section~\ref{notations} for the precise normalization).

\subsubsection{Oscillatory integrals}

Our aim in this paper is to confront pseudo-products with  oscillatory integrals, about which we now say a word.

\mb The most simple instance is certainly that of real-valued linear forms $f \mapsto \int f(x) m(x) e^{i\lambda \phi(x)} \,dx$. If the Hessian of $\phi$ is nowhere degenerate, the behaviour of this expression as $\lambda$ goes to infinity is described by the stationary phase lemma; the other possibilities are more subtle.

\mb The next step is given by linear oscillatory integrals maps between function spaces. The prototypes are $Lf(x) = \int e^{i\lambda \phi(x,y)} f(y) m(x,y)\,dy$ (oscillations in physical space) and $Lf(x) = \int e^{ix\cdot\xi} e^{i\lambda P(\xi)} \widehat{f}(\xi) \,d\xi$ (oscillations in Fourier space). This last operator corresponds to solutions of the linear equation $i \partial_t f + P(D) f = 0$, see the next section. We refer to the books of Stein~\cite{Stein} and Sogge~\cite{Sogge} for a discussion of the above operators, and other instances where linear oscillatory integrals occur, such as the theory of Fourier integral operators, or the theory of Fourier restriction.

\mb 
Finally, multilinear oscillatory integrals of the type
$$
(f_1,\dots,f_n) \mapsto  \int e^{i\lambda \phi(x_1,\dots,x_n)} f_1(x_1) \dots f_n(x_n) \,dx_1\dots dx_n
$$
were recently considered by several authors, we mention in particular Phong, Stein and Sturm~\cite{PSS} and Christ, Li, Tao and Thiele~\cite{CLTT}.

\subsubsection{Bringing them together}

We want to consider here the following instance of multilinear oscillatory integral operators: pseudo-product operators whose symbol contains an oscillatory phase. We shall simply consider the bilinear case, where the operator reads
\be{def}
B_\lambda (f,g) \overset{def}{=} \mathcal{F}^{-1} \int_{\mathbb{R}^d} m(\eta,\xi-\eta) e^{i\lambda \phi(\eta,\xi-\eta)} \widehat{f}(\eta) \widehat{g} (\xi-\eta) \,d\eta.
\ee
The question that we ask is the following : under which conditions on $m$ and $\phi$ is $B_\lambda$ bounded between Lebesgue spaces ? How does the bound depend on $\lambda$ ?

\subsection{Application to PDEs}

\subsubsection{Long term behaviour for a general dispersive PDE}
Consider a general nonlinear quadratic dispersive PDE
$$ 
\left\{ \begin{array}{l}
\partial_t u + iP(D) u = T_m(u,u) \\
u(t=0)=u_0.
\end{array} \right.
$$
where we follow the above notation in denoting $T_m$ for the pseudo-product with symbol $m$; of course, nonlinearities of higher order can be dealt with in a similar way to what we will explain. We do not consider such a nonlinearity for the sake of generality, but because it does actually occur in PDE problem. For instance, the nonlinearity of the water-waves problem can be expanded as a sum of pseudo-product operators: see~\cite{GMS2}. Or it is well-known that a nonlinearity $H(u)$ can for many purposes be replaced by its paralinearization $ H(u(t,.)) \simeq \pi_{H'(u(t,.))} (u(t,.))$: see the seminal work of Bony~\cite{bony}).

Let us now write Duhamel's formula for the solution of the above problem: it reads
$$
u(t) = e^{itP(D)} u_0 + \int_0^t e^{i(s-t)P(D)} T_m(u(s),u(s))\,ds.
$$
Our aim is to understand how $u$ behaves for large $t$, in particular whether it scatters.

\subsubsection{The linear part}
The linear part of the above right-hand side satisfies the dispersive estimates alluded to above. For instance, if $P(D)=\Delta$, then, for $p\in[1,2]$,
\be{strichartz} \left\| e^{it\Delta}(f) \right\|_{p'} \lesssim |t|^{-\frac{d}{2}\left(\frac{1}{p}-\frac{1}{p'}\right)} \|f\|_p. \ee
(On the one hand the $p=2$ inequality is a direct consequence of Plancherel's equality and on the other hand the $p=1$ inequality is a consequence of the stationary phase Lemma. The intermediate exponents $p$ are deduced by interpolation.) 

These dispersive estimates yield in turn Strichartz estimates: see Ginibre and Velo~\cite{GV}, and Keel and Tao~\cite{KT}.

\subsubsection{The bilinear part}
We are interested in reproducing a similar reasoning for the bilinear term, namely we want to understand when it is bounded in various space time norms, and in particular how it decays as $t$ goes to infinity. The most simple possibility is to use the boundedness of $T_m$ in appropriate spaces, and the linear estimates, but it only gives a partial answer. It is indeed possible to obtain sharper results if one is willing to work in a more authentically bilinear way: first instances of this approach go back at least to the normal form method of Shatah~\cite{Shatah} and the vector field method of Klainerman~\cite{Klainerman}. Following linear Strichartz estimates, bilinear Strichartz estimates have been developped, and proved very useful: see for instance Klainerman and Foschi~\cite{KF}.

In order to better understand the bilinear term in the above Duhamel equaion, let us change the unknown function from $u$ to 
$$ f(t,x) = e^{itP(D)}[u(t,.)](x).$$
The idea is the following: in the regime where the equation scatters (which we focus on), $f$ is converging as $t$ goes to infinity, whereas $u$ is not, due to oscillations in frequency space. Thus, by writing the Duhamel term
\begin{equation}
\label{duhamel}
\int_0^t e^{i(s-t)P(D)} T_m(u(s),u(s))\,ds = \int_0^t \int \int e^{ix\cdot(\xi+\eta)} e^{i[(s-t)P(\xi+\eta) - s P(\eta) - s P(\xi)]} \widehat{f}(\eta,s) \widehat{f}(\xi,s) \,d\xi\,d\eta\,ds,
\end{equation}
we isolate in the right-hand side the oscillations in the term $e^{i[(t-s)P(\xi+\eta) - s P(\eta) - s P(\xi)]}$. 

The relation to $B_\lambda$ is now clear, and why we believe its boundedness properties can have implications for the PDE theory: our work essentially enables one to understand the behaviour of the integrand (in $s$) of the right-hand side of~(\ref{duhamel}).

What is missing is also clear: understanding the effects of the $s$ integration. The implications are best understood in physical terms if one uses the notions of \textit{space resonance} and \textit{time resonance}: by ignoring the $s$ integration, one ignores the effects of time resonances, and focuses on space resonances. For an explanation of these concepts, as well as applications of these to nonlinear PDE, we refer the reader to works of the second author, Masmoudi and  Shatah \cite{GMS1,GMS2,GMS3}.

\subsubsection{A restriction on the phase function} It is important to notice that the phase function occuring in a PDE setting (such as above) is much less general than what was considered above: it is a sum of three functions, respectively of $\xi$, $\eta$, and $\xi+\eta$. For such a particular phase functions, some of the theorems which follow can be proved with considerably less effort.

\subsection{Results}

We only treat the bilinear case. However we emphasize that even in the linear case, the estimates involving a
mixture between an oscillatory term and a Coifman-Meyer symbol seem to be new.

\mb
We begin by describing in \textit{Section \ref{sec:1}} the most simple case: polynomial phase $\phi$ of order $2$ and $d=1$. Then we give in \textit{Section \ref{sec:2}} boundedness results for (\ref{def}) for smooth phase and symbol in Lebesgue spaces and weighted Lebesgue spaces. More precisely, we obtain two kind of estimates, first
$$ \left\|B_\lambda(f,g)\right\|_{L^\infty} \lesssim |\lambda|^{-d} \|f\|_{L^1} \|g\|_{L^1}$$
and then
$$ \left\|B_\lambda(f,g)\right\|_{L^2} \lesssim |\lambda|^{-d/2} \|f\|_{L^2} \|g\|_{L^1}.$$
We obtain a full set of inequalities by interpolating between these two estimates in Theorem \ref{thm:thm3smooth} ( and a weighted version of these results in Theorem \ref{thm:thm1smoothw}). \\
Similar results for smooth symbols $m$ supported on a submanifold are developped in \textit{Section \ref{sec:singular}}. This case seems to make appear some very difficult questions. \\
In \textit{Section~\ref{sec:3}}, we are specially interested in proving similar estimates for ($x$ independent) Coifman-Meyer symbols $m$: this is achieved in Theorem \ref{thm:thm3cm}.
These estimates are then extended to $x$-dependent non-smooth symbols in \textit{Section \ref{sec:4}}. \\
\textit{Section \ref{sec:5}} is devoted to the proof of optimality for our estimates in the following sense: the range of exponents obtained by interpolation between the $L^1\times L^1\to L^\infty$ and $L^2\times L^1\to L^2$ estimates is the biggest one where boundedness can be obtained. We finish our work by describing in \textit{Section \ref{sec:6}} an application of these bilinear oscillatory integrals in order to prove boundedness of some bilinear multipliers (associated to non-smooth symbols) in products of Lebesgue spaces.

\section{Notations} \label{notations}

We sometimes denote $C$ for a constant whose value may change from one line to the other. Mostly however, we use $\gtrsim$ and $\lesssim$: given two quantities $A$ and $B$, we write $A \gtrsim B$ if there exists a constant $C$ such that $A \geq C B$; there is an obvious adaptation for $\lesssim$.

We denote $A >> B$ if there exists a big enough (depending on the context) constant $C$ such that $A \geq CB$.

Given a (real) function $\phi(\xi,\eta)$, with $(\eta,\xi) \in \mathbb{R}^d\times \mathbb{R}^d$, we denote $\nabla_\xi \nabla_\eta \phi$ for the matrix
$$
\nabla_\xi \nabla_\eta \phi \overset{def}{=} \left( \frac{\partial^2 \phi}{\partial \xi_i \partial \eta_j} \right)_{i,j}.
$$
The Hessian of $\phi$ can then be written as
$$
\operatorname{Hess} \phi = \left( \begin{array}{ll} \nabla_\xi \nabla_\xi \phi & \nabla_\xi \nabla_\eta \phi \\ \nabla_\xi \nabla_\eta \phi & \nabla_\eta \nabla_\eta \phi \\ \end{array} \right).
$$
The Fourier transform of $f$ is denoted by $\mathcal{F} f$ or $\widehat{f}$, and defined as follows
$$
\mathcal{F} f (\xi) = \frac{1}{(2\pi)^{d/2}} \int e^{- i x \cdot \xi} f(x) dx \,\,.
$$
For $1 \leq p \leq \infty$, the $L^p$ norm of $f$ is denoted by $\|f\|_{L^p(\mathbb{R}^d)}$, or simply $\|f\|_p$, and defined by
$$
\|f\|_p = \left[ \int_{\mathbb{R}^d} |f(x)|^p \,dx \right]^{1/p}
$$
with the usual modification if $p = \infty$. For $p\in[1,\infty]$ and $a\geq 0$, we denote the weighted Lebesgue space $L^p(\langle x \rangle ^a)$ corresponding to the norm
$$ \|f\|_{L^p(\< x \> ^a)} \overset{def}{=} \|x\to \< x \> ^a f\|_{L^p},$$
with as usual $\<x\> \overset{def}{=}(1+|x|^2)^{1/2}$. \\
We denote by $\mathcal{H}^1$ for the classical Hardy space on $\R^{d}$ (see the initial work of R. Coifman and G. Weiss \cite{CW} for its first definition and the book of E. Stein \cite{Stein} for the study and several characterizations).

\begin{df} Let $T$ be a bilinear operator bounded from $L^{p_1} \times L^{p_2}$ into $L^p$ for exponents $p,p_1,p_2\in[1,\infty]$. Using real duality, we define its two adjoints $T^{*1}$ and $T^{*2}$ by
$$ \langle T(f,g) , h \rangle \overset{def}{=} \langle T^{*1}(h,g),f\rangle \overset{def}{=} \langle T^{*2}(f,h),g\rangle.$$
So $T^{*1}$ is bounded from $L^{p'} \times L^{p_2}$ into $L^{p_1'}$ and $T^{*2}$ is bounded from $L^{p_1} \times L^{p'}$ into $L^{p_2'}$.
\end{df}

\section{The simple case : $\phi$ is a polynomial function of order $2$ and $d=1$} \label{sec:1}

Let us treat in this section the particular case of a polynomial phase $\phi$ of degree less than $2$ and one dimensional variables.
Using the modulation invariance of the Lebesgue norms, we have only to deal with the homogeneous polynomial phase $\phi$ of degree $2$. So we are working with real variables and the phase $\phi$ takes the form:
$$ \phi(\eta,\xi)\overset{def}{=}a\eta^2+b\eta\xi+c\xi^2.$$
We write it in the following canonical form, involving only $\eta^2$, $\xi^2$ and $(\eta+\xi)^2$:
\be{factorphi} \phi(\eta,\xi)=\frac{b}{2}(\eta+\xi)^2 + \left(a-\frac{b}{2}\right) \eta^2 + \left(c-\frac{b}{2}\right) \xi^2. \ee

\mb Let us work with a bilinear multiplier $T_m$, given by a bilinear symbol $m$, belonging to the space ${\mathcal M}_{p,q,r}$ in the local-$L^2$ case. This means that for all exponents $p,q,r$ satisfying the homogeneous relation
$$\frac{1}{p}+\frac{1}{q}=\frac{1}{r}$$
and the local-$L^2$ condition: $2\leq p,q,r'<\infty$, the bilinear operator $T_m$ is bounded from $L^p(\R) \times L^q(\R)$ into $L^r(\R)$.

\begin{ex} Several classes of bilinear multipliers satisfy to this property~:
\begin{itemize}
\item The paraproducts and the Coifman-Meyer multipliers (see the works of Bony in \cite{bony} and of  Coifman and Meyer in \cite{CM2,coifmanmeyer,CM3} and \cite{B0} for some uniform estimates and Grafakos and Torres~\cite{GTo} for boundedness with exponents in the whole optimal range).
\item The operators with flag singularities (Muscalu~\cite{M}).
\item The Marcinkiewicz multipliers under some assumptions (see the work of Grafakos and Kalton \cite{GK}).
\item The multiparameter paraproducts (Muscalu, Pipher, Tao and Thiele~\cite{MPTT}).
\item The bilinear Hilbert transforms and related bilinear multipliers with modulation symmetry (see the works \cite{LT2,LT1,LT3,LT4} of Lacey and Thiele, \cite{GN1,GN2,GN3} of Gilbert and Nahmod, and \cite{mtt,MTT2} of Muscalu,  Tao and Thiele).
 \item The indicator function of the unit disc (see the work \cite{GL} of Grafakos and Li).
\end{itemize}
\end{ex}

\begin{thm} \label{thm:casfacile} Let us assume that
\be{condition1} b\overset{def}{=}\partial_\eta \partial_\xi \phi \neq 0, \quad 2a-b = \left(\partial_\eta^2-\partial_\eta \partial_\xi\right) \phi \neq 0 \textrm{  and  } 2c-b = \left(\partial_\xi^2-\partial_\eta \partial_\xi\right) \phi \neq 0. \ee
Then the bilinear oscillatory integral
$$B_\lambda(f,g)(x)\overset{def}{=} \frac{1}{(2\pi)^{d/2}}\int_{\R^2} e^{ix\cdot (\eta+\xi)} e^{-i\lambda \phi(\eta,\xi)} m(\eta,\xi) \widehat{f}(\eta)\widehat{g}(\xi) d\eta d\xi$$ satifies to the following boundedness~: for all exponents $p,q,r' \in (1,2]$ verifying
$$\frac{1}{r}+1=\frac{1}{p}+\frac{1}{q},$$ there exists a constant $C=C(p,q,r,\phi,m)$ such that for all $\lambda \neq 0$
$$ \left\|B_\lambda(f,g) \right\|_{L^r} \leq C|\lambda|^{-\frac{1}{2}} \|f\|_{L^p} \|g\|_{L^q}.$$
\end{thm}

\dem According to (\ref{factorphi}), the oscillatory integral $B_\lambda$ can be written as follows~:
$$B_\lambda(f,g) \overset{def}{=} (2\pi)^{-d/2} e^{-it\frac{b}{2} \Delta} T_{m} \left( e^{-it\left(a-\frac{b}{2}\right)\Delta} f , e^{-it\left(c-\frac{b}{2}\right)\Delta} g \right).$$
Then the results are a direct consequence of the classical dispersive estimates (\ref{strichartz}) and of the boundedness of $T_{m}$ in the local-$L^2$ case from $L^{p'} \times L^{q'}$ to $L^{r'}$. \findem

\mb We leave to the reader the corresponding results when in (\ref{condition1}) only one or two terms are vanishing. Moreover if we know some boundedness of $T_{m}$ with some infinite exponents, then we can allow $p=1$ or $q=1$ in Theorem \ref{thm:casfacile}.

\section{The case of a smooth phase and symbol} \label{sec:2}

This section is devoted to the particular case where both the phase $\phi$ and the symbol $m$ are supposed to be smooth and compactly supported on $\R^{2d}$.

\mb First we deal with estimates in classical Lebesgue spaces. Then we study the behavior of the oscillatory integrals in weighted Lebesgue spaces.

\subsection{Estimates on Lebesgue spaces}

We have different kinds of estimates for the considered bilinear oscillatory integral~:
$$ B_\lambda(f,g)(x) \overset{def}{=} \frac{1}{(2\pi)^{d/2}}\int_{\R^{2d}} e^{ix\cdot (\eta+\xi)} e^{i\lambda\phi(\eta,\xi)} m(\eta,\xi) \widehat{f}(\eta) \widehat{g}(\xi) d\eta d\xi.$$
The first one describes a $L^1 \times L^1\to L^\infty$ decay~:

\begin{thm} \label{thm:thm1smooth} Consider $\phi,m\in C^\infty_0(\R^{2d})$ and assume that the Hessian matrix $\operatorname{Hess}(\phi)$ is non degenerate on $\operatorname{Supp}(m)$. Then there exists an implicit constant such that for all $\lambda\neq 0$~:
$$ \forall f,g\in \s(\R^d), \qquad \left\|B(f,g)\right\|_{L^\infty} \lesssim |\lambda|^{-d} \|f\|_{L^1} \|g\|_{L^1}.$$
\end{thm}

\dem The bilinear oscillatory integral $B_\lambda$ can be seen as a bilinear operator whose kernel $K$ reads
\be{noyau} K_\lambda (x,y,z) \overset{def}{=} \frac{1}{(2\pi)^{d}} \int_{\R^{2d}} e^{i\eta.(x-y)} e^{i\xi.(x-z)} m(\eta,\xi) e^{i\lambda\phi(\eta,\xi)} d\eta d\xi.\ee
Let us denote the new phase
$$\tilde{\phi}(\eta,\xi)\overset{def}{=}\phi(\eta,\xi) + \lambda^{-1}\left[\eta\cdot (x-y)+\xi\cdot (x-z)\right].$$
The assumption yields that the Hessian matrix $\operatorname{Hess}(\tilde{\phi})=\operatorname{Hess}(\phi)$ is non degenerate on $\operatorname{Supp}(m)$. The stationary phase lemma gives
$$ \left|K(x,y,z) \right| = \frac{1}{(2\pi)^{d}} \left|\int_{\R^{2d}} e^{i\lambda \tilde{\phi}(\eta,\xi)} m(\eta,\xi) d\eta d\xi \right| \lesssim |\lambda|^{-d}
$$
with an implicit constant dependent only on $\operatorname{Hess}(\phi)$. Thus $\|K_\lambda\|_{L^\infty}\lesssim \lambda^{-d}$ which gives the desired result. \findem

\mb To deal with $\phi$ which are non degenerate in one direction only, one can use Proposition 5 of Chapter VIII in Stein~\cite{Stein} instead of the classical stationary phase Lemma, to obtain the following extension~:

\begin{prop} \label{prop:prop1smooth} Consider $\phi,m\in C^\infty_0(\R^{2d})$ and assume that for a multi-index $\alpha\in (\N^{2d})^k$ with $k=|\alpha|\geq 2$, we have
 $$\left|\partial^\alpha \phi(\eta,\xi) \right|  \gtrsim 1$$ on $\operatorname{Supp}(m)$. Then there exists an implicit constant such that for all $\lambda\neq 0$~:
$$ \forall f,g\in \s(\R^d), \qquad \left\|B(f,g)\right\|_{L^\infty} \lesssim |\lambda|^{-\frac{1}{2|\alpha|}} \|f\|_{L^1} \|g\|_{L^1}.$$
\end{prop}

\dem The proof is exactly the same as the previous one, in using Proposition 5 of Chapter VIII in Stein~\cite{Stein}) instead of the stationary phase Lemma and the fact that for all multi-index $\beta$ with $|\beta|=|\alpha|$, we have~:
$$ \partial^\beta \tilde{\phi}=\partial^\beta \phi.$$
We remark that in Proposition 5 of Chapter VIII in~\cite{Stein}, the implicit constant, written $c_k(\phi)$ is in fact bounded by the homogeneous norm $\|\nabla^{|\alpha|}(\phi)\|_\infty$ and not only by the inhomogeneous norm $\|\phi\|_{C^{|\alpha|}}$. That is why we can apply this result. \findem

\begin{thm} \label{thm:thm2smooth} Consider $\phi,m\in C^\infty_0(\R^{2d})$ and assume that the derivatives-matrix $\nabla_\xi (\nabla_\eta-\nabla_\xi) \phi$ is non degenerate on $\operatorname{Supp}(m)$. Then there exists an implicit constant such that for all $\lambda\neq 0$~:
$$ \forall f,g\in \s(\R^d), \qquad \left\|B_\lambda(f,g)\right\|_{L^2} \lesssim |\lambda|^{-d/2} \|f\|_{L^2} \|g\|_{L^1}.$$
 \end{thm}

\dem Computing the Fourier transform of $B$, we get
\begin{align*}
 \widehat{B_\lambda(f,g)}(\xi) & \overset{def}{=}\int_{\R^d} e^{i\lambda \phi(\eta,\xi-\eta)} m(\eta,\xi-\eta) \widehat{f}(\eta) \widehat{g}(\xi-\eta) d\eta \\
 & = \frac{1}{(2\pi)^{d/2}} \int_{\R^d} g(x) \left[\int_{\R^d} e^{-ix\cdot (\xi-\eta)} e^{i\lambda \phi(\eta,\xi-\eta)} m(\eta,\xi-\eta) \widehat{f}(\eta)  d\eta \right] dx.
\end{align*}
So let us denote by $T_\lambda^x$ the linear operator
$$T_\lambda^x(h)(\xi) \overset{def}{=} \int_{\R^d} e^{-ix\cdot (\xi-\eta)} e^{i\lambda \phi(\eta,\xi-\eta)} m(\eta,\xi-\eta) h(\eta)  d\eta.$$
By appealing to Plancherel's Theorem, it suffices to prove
\be{eq:thm1smooth} \sup_{x\in\R^d} \left\|T_\lambda^x(h)\right\|_{L^2} \lesssim |\lambda|^{-d/2}\|h\|_{L^2}. \ee
To estimate in $L^2$ the operator $T_x$, we refer the reader to Proposition 1.1 of Chapter IX in Stein~\cite{Stein} for a detailed proof. Moreover we refer the reader to Theorem \ref{thm:thm2cm}, whose proof is detailed and contains all the arguments, though in a more complex framework. For an easy reference, we quickly recall the ideas. \\
By a $TT^*$ argument, it suffices to prove the bound $\sup_{x\in\R^d} \left\|[T_\lambda^x (T_\lambda^x)^*](h)\right\|_{L^2} \lesssim |\lambda|^{-d/2}\|h\|_{L^2}4$. Then, using integrations by parts and bounding the resulting expressions by the assumption on $\nabla_\xi (\nabla_\eta - \nabla_\xi) \phi$, it can be proved that the kernel $K$ of the operator $T^x(T^x)^*$ satisfies
\be{eq:ker} \left|K(\xi,\eta)\right| = \left|\int_{\R^d} \overline{m(\tau,\eta-\tau)} m(\tau,\xi-\tau) e^{i\lambda(\phi(\tau,\xi-\tau)-\phi(\tau,\eta-\tau))} d\tau \right| \lesssim \left(1+\lambda|\xi-\eta|\right)^{-N} \ee
for every large enough integer $N$. The desired bound follows. \findem

\mb Having obtained two kinds of bilinear estimates, we can interpolate between them. A simple computation gives that
$$ B_\lambda^{*1}(h,g)(x)\overset{def}{=}\frac{1}{(2\pi)^{d/2}} \int_{\R^{2d}} e^{ix\cdot (\eta+\xi)}e^{i\phi(\eta-\xi,\xi)} \widehat{h}(\eta) \widehat{g}(\xi) m(\eta-\xi,\xi) d\eta d\xi.$$
So it corresponds to the bilinear oscillatory integral associated to the phase and symbol
$$ \phi^{*1}(\eta,\xi)\overset{def}{=}\phi(\eta-\xi,\xi)\quad \textrm{ and } \quad m^{*1}(\eta,\xi)\overset{def}{=}m(\eta-\xi,\xi).$$
Similarly $B_\lambda^{*2}$ corresponds to the bilinear oscillatory integral associated to the phase and symbol
$$ \phi^{*2}(\eta,\xi)\overset{def}{=}\phi(\eta,\xi-\eta)\quad \textrm{ and } \quad m^{*2}(\eta,\xi)\overset{def}{=}m(\eta,\xi-\eta).$$

\begin{rem} \label{rem} We remark that $\operatorname{Hess}(\phi)(\eta,\xi)$ is non degenerate if and only if $\operatorname{Hess} (\phi^{*1})(\eta+\xi,\xi)$ is non degenerate if and only if $\operatorname{Hess}(\phi^{*2})(\eta,\xi+\eta)$ is non degenerate (indeed the determinant of the three Hessian matrices are equal). Moreover we have
$$ \nabla_\xi (\nabla_\eta-\nabla_\xi) \phi^{*1}(\eta,\xi) = \left[ \left(2\nabla_\eta - \nabla_\xi\right)(\nabla_\xi-\nabla_\eta) \phi \right] (\eta-\xi,\xi) $$
and
$$\nabla_\xi (\nabla_\eta-\nabla_\xi) \phi^{*2}(\eta,\xi) = \left[ \left(\nabla_\eta - 2\nabla_\xi\right)\nabla_\xi \phi \right] (\eta,\xi-\eta).$$
\end{rem}

\begin{thm} \label{thm:thm3smooth} Consider $\phi,m\in C^\infty_0(\R^{2d})$ and assume that $\operatorname{Hess}(\phi)$, $\nabla_\xi (\nabla_\eta-\nabla_\xi) \phi$, $\left(2\nabla_\eta - \nabla_\xi\right)(\nabla_\xi-\nabla_\eta) \phi$ and $\left(\nabla_\eta - 2\nabla_\xi\right)\nabla_\xi \phi$ are non degenerate on $\operatorname{Supp}(m)$. \\
Then for all exponents $p,q,r$ verifying
\begin{equation} \label{eq:condi}
\left\{ \begin{aligned}
& \frac{1}{p} + \frac{1}{q} + \frac{1}{r} \leq 2 \\
& \frac{1}{p} + \frac{1}{q} - \frac{1}{r} \geq 1 \\
& \frac{1}{p} - \frac{1}{q} - \frac{1}{r} \leq 0 \\
& \frac{1}{p} - \frac{1}{q} + \frac{1}{r} \geq 0 \\
\end{aligned} \right.
\end{equation}
there exists a constant $C=C(p,q,r,\phi,m)$ such that for all $\lambda\neq 0$
$$ \left\|B_\lambda(f,g) \right\|_{L^r} \leq C|\lambda|^{-\frac{d}{2}\left(\frac{1}{p}+\frac{1}{q}-\frac{1}{r}\right)} \|f\|_{L^p} \|g\|_{L^q}.$$
\end{thm}

The set of $(\frac{1}{p},\frac{1}{q},\frac{1}{r})$ satisfying the above inequalities is not symmetrical; but the set given by the triplets $(\frac{1}{p},\frac{1}{q},\frac{1}{r'})$ (which corresponds to considering the trilinear form associated by duality) such that $(p,q,r)$ is admissible is symmetrical. Therefore we choose to represent it below.

\begin{figure}[htbp]
 \centering
\includegraphics[width=0.5\textwidth]{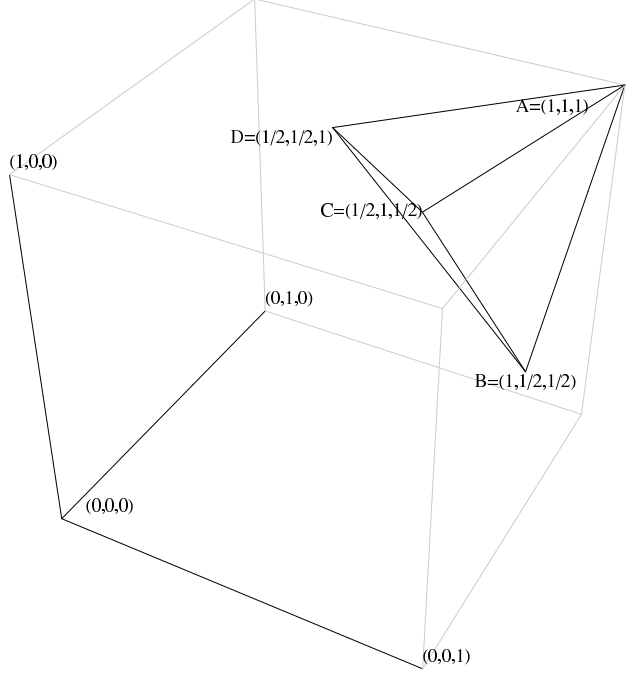}
\caption{The coordinates are $(\frac{1}{p},\frac{1}{q},\frac{1}{r'})$; the solid tetrahedron corresponds to the $(p,q,r)$ which satify the inequalities in the theorem.} \label{fig:tetra}
\end{figure}

\dem Consider the trilinear form
$$ T(f,g,h)\overset{def}{=} \langle B_\lambda(f,g),h\rangle = \langle f,B_\lambda^{*1}(h,g) \rangle = \langle g,B_\lambda^{*2}(f,h)\rangle.$$
The assumptions and Remark \ref{rem} permit to apply Theorems \ref{thm:thm1smooth} and \ref{thm:thm2smooth} to the operators $B,B^{*1}$ and $B^{*2}$. So by duality, we deduce the following boundedness for $T$~:
\begin{align*}
 \left|T(f,g,h)\right| & \lesssim |\lambda|^{-d} \|f\|_{L^1}\|g\|_{L^1} \|h\|_{L^1} \\
 \left|T(f,g,h)\right| & \lesssim |\lambda|^{-d/2} \|f\|_{L^1}\|g\|_{L^2}\|h\|_{L^2} \\
 \left|T(f,g,h)\right| & \lesssim |\lambda|^{-d/2} \|f\|_{L^2}\|g\|_{L^1}\|h\|_{L^2} \\
 \left|T(f,g,h)\right| & \lesssim |\lambda|^{-d/2} \|f\|_{L^2}\|g\|_{L^2}\|h\|_{L^1}.
\end{align*}
Then we can now use trilinear interpolation and deduce the desired estimates. We refer the reader to \cite{BL} for a multilinear version of the Riesz-Thorin theorem (complex interpolation) and to \cite{GM,GT,mtt} for a multilinear theory of real interpolation. 
\findem

\subsection{Estimates on weighted Lebesgue spaces}

In this section, we are looking for a weighted version of the previous results.

\begin{thm} \label{thm:thm1smoothw} Consider $\phi,m\in C^\infty_0(\R^{2d})$ $a,b\geq 0$ and assume that the Hessian $\operatorname{Hess} (\phi)$ and $(\nabla_\eta-\nabla_\xi)\phi$ are non degenerate on $\operatorname{Supp}(m)$. Then there exists an implicit constant such that for all $|\lambda|>1 0$~:
$$ \forall f,g\in \s(\R^d), \qquad \left\|B_\lambda(f,g)\right\|_{L^\infty(\< x\> ^{a})} \lesssim |\lambda|^{-d-b} \|f\|_{L^1( \< x\> ^{a+b})} \|g\|_{L^1(\< x\> ^{a+b})}.$$
\end{thm}

\dem We follow the proof of Theorem \ref{thm:thm1smooth}. \\
$\bullet$ If $a=0$ and $b\in\N$, then we have to consider (instead of (\ref{noyau})) the new kernel given by
$$ K(x,y,z) = \< y \>^{-b}\< z \> ^{-b}\int_{\R^{2d}} e^{i\eta\cdot (x-y)} e^{i\xi\cdot (x-z)} m(\eta,\xi) e^{i\lambda\phi(\eta,\xi)} d\eta d\xi .$$
Changing the integration variables to $\alpha=\xi+\eta$ and $\beta=\eta-\xi$, one gets:
$$ \left|K(x,y,z)\right| \lesssim \< y \>^{-b}\< z \> ^{-c}\int_{\R^{2d}} e^{i\alpha\cdot x} e^{i(\alpha+\beta)\cdot \frac{y}{2}}e^{i(\alpha-\beta)\cdot \frac{z}{2}}  m(\frac{\alpha+\beta}{2},\frac{\alpha-\beta}{2}) e^{i\lambda\phi(\frac{\alpha+\beta}{2},\frac{\alpha-\beta}{2})} d\alpha d\beta.$$
Use now the identity
$$
\frac{\nabla_\beta \phi}{2 \lambda |\nabla_\beta \phi|^2} \cdot \nabla_\beta e^{i \lambda \phi(\frac{\alpha+\beta}{2},\frac{\alpha-\beta}{2})} = e^{i \lambda \phi(\frac{\alpha+\beta}{2},\frac{\alpha-\beta}{2})}
$$
to integrate by parts $b$ times with respect to $\beta$; and use the stationary phase lemma to estimate the resulting expression. The outcome is the bound
\begin{align*}
 \left|K(x,y,z)\right| & \lesssim \< y \> ^{-b} \< z \> ^{-b}|\lambda|^{-d-b} (\< y \> + \< z\>) ^b \\
 & \lesssim |\lambda|^{-d-b},
\end{align*}
which is the desired estimate if $a=0$, and $b$ is an integer. \\
$\bullet$ By interpolation between weighted Lebesgue spaces (see \cite{SteinWeiss}), this result is extended to positive real $b$. \\
$\bullet$ Suppose now that $a$ is an integer. Observe that it suffices to prove the estimate of the theorem with the weight $<x>^a$ on the right-hand side replaced by $|x^A|= |x_1^{a_1} x_2^{a_2} \dots x_n^{a_n} |$, for a multiindex $A=(a_1 \dots a_n)$ of size $|A| = \alpha_1 + \dots + \alpha_n$ less than $a$. \\
Multiplying $B_\lambda(f,g)$ by $x^A$ corresponds, in Fourier space, to applying $\partial_\xi^\alpha = \partial_{\xi_1}^{\alpha_1} \dots \partial_{\xi_n}^{\alpha_n}$ to $\int m(\eta,\xi-\eta) e^{i\lambda \phi(\eta,\xi-\eta)} \widehat{f}(\eta) \widehat{g}(\xi-\eta)\,d\eta$. It is harmless if the derivatives hit $m$, thus we consider that they always hit either the oscillating factor, or $\widehat{g}(\xi-\eta)\,d\eta$. Denoting indifferently $\partial_\xi$ for any partial derivative in $\xi$, it means that
$\mathcal{F} \left[x^A B_\lambda(f,g) \right]$, up to easily bounded terms, reads
$$
\sum_{\ell=1}^{|A|} \int m(\eta,\xi-\eta)\left[\lambda \partial_\xi \phi(\eta,\xi-\eta)\right]^\ell e^{i\lambda \phi(\eta,\xi-\eta)} \widehat{f}(\eta) \partial_\xi^{|A|-\ell} \widehat{g}(\xi-\eta)\,d\eta.
$$
The Fourier transform of this sum can be rewritten
$$
\sum_{\ell=1}^{|A|} \lambda^\ell B_\lambda(f,x^{|A|-\ell} g)
$$
where the different $B_\lambda$ occuring in the above sum have different (smooth) symbols. By the case $a=0$ discussed aboved, this sum can easily be bounded in $L^\infty$:
\begin{equation}
\begin{split}
\left\|\sum_{\ell=1}^{|A|} \lambda^\ell B_\lambda(f,x^{|A|-\ell} g) \right\|_{L^\infty} & \lesssim \sum_{\ell=1}^{|A|} \lambda^\ell \lambda^{-d-b-\ell} \|f\|_{L^1( \< x\> ^{b+\ell})} \|x^{|A|-\ell} g\|_{L^1(\< x\> ^{b+\ell})} \\
& \lesssim \lambda^{-d-b} \|f\|_{L^1( \< x\> ^{a+b})} \|g\|_{L^1( \< x\> ^{a+b})},
\end{split}
\end{equation}
which is the desired bound. \\
$\bullet$ We conclude for non-negative real $a$ by interpolating again. \findem

\mb For the $L^1\times L^2 \to L^2$ decay, we get the following weighted version~:

\begin{prop} \label{prop:prop2smoothw} Consider $\phi,m\in C^\infty_0(\R^{2d})$, $a,b\geq 0$ and assume that the derivative-matrix $\nabla_\xi(\nabla_\eta-\nabla_\xi) \phi$ and $(\nabla_{\eta}-\nabla_{\xi})\phi$ are non degenerate on $\operatorname{Supp}(m)$. Then there exists an implicit constant such that for all $\lambda\neq 0$~:
$$ \forall f,g\in \s(\R^d), \qquad \left\|B_\lambda(f,g)\right\|_{L^2(\< x\> ^{a})} \lesssim |\lambda|^{-d/2-b} \|f\|_{L^2(\< x\> ^{a+b})} \|g\|_{L^1(\< x\> ^{a+b})}.$$
\end{prop}

\mb We let the proof to the reader in combining the proof of Theorem \ref{thm:thm2smooth} and the weighted proof for Proposition \ref{prop:prop1smooth}.

\mb As previously, by interpolating the weighted Lebesgue spaces, we can obtain boundedness for other exponents with some appropriate decay in $|\lambda|$.

\begin{rem} The above weighted estimate has a very natural interpretation in terms of partial differential equations, which we now explain; it corresponds to the space resonance phenomenon observed in~\cite{GMS1} and then used in~\cite{GMS2}~\cite{GMS3} in order to understand the long term interactions between waves. \\
Consider $u$ and $v$, solutions of a linear dispersive equation (with dispersion relation $\tau = P(\xi)$)
$$
\left\{ \begin{array}{l} i \partial_t u + P(D) u = 0 \\ u(t=0) = u_0 \end{array} \right. \quad \left\{ \begin{array}{l} i \partial_t v + P(D) v = 0 \\ v(t=0) = v_0 \end{array} \right.
$$
and assume that $u_0$, and $v_0$ are localized in physical space around $0$, and in frequency around respectively $\xi_0$ and $\xi'_0$. It is well known that the resulting wave packets $u$ and $v$ will be localized respectively around $x \sim - \nabla P(\xi_0) t$, and $x \sim - \nabla P(\xi'_0)t$. If one considers the product $u(t)v(t)$ (of course, things would be nearly identical for a general pseudo-product operator), it will be small if $P'(\xi_0) \neq P'(xi'_0)$. In physical terms: wave packets with different group velocities do not interact much. \\
To see how this is linked to the above propositions, observe that $u(t)v(t)$ can be written
$u(t)v(t) = B_t(u_0,v_0)$ where $B$ has symbol identically equal to $1$, and phase $\phi(\xi,\eta) = P(\eta) + P(\xi)$. The condition that $(\nabla_\eta - \nabla_\xi) \phi$ essentially means $P'(\xi_0) \neq P'(\xi'_0)$, or in other words: group velocities do not coincide. Then, the above propositions can be read as ``localization of the data leads to better decay'', which is the quantitative version of the physical fact explained above.
\end{rem}

\section{The case of symbols supported on a submanifold} \label{sec:singular}

In this section, we are studying bilinear oscillatory integrals associated to symbols $m$ supported on a manifold of $\R^{2d}$. So let $\Gamma$ a smooth manifold of $\R^{2d}$ and consider $\sigma_\Gamma$ its superficial measure and denote $\delta=dim(\Gamma)$ its dimension. \\
We consider the following bilinear integral~:
$$B_\lambda(f,g)(x)\overset{def}{=}\frac{1}{(2\pi)^{d/2}} \int_\Gamma e^{ix\cdot (\eta+\xi)} \widehat{f}(\eta) \widehat{g}(\xi) m(\eta,\xi) e^{i\lambda\phi(\eta,\xi)}
d\sigma_\Gamma(\eta,\xi).$$

\begin{thm} \label{thm:thm1curve} Consider $\phi\in C^\infty_0(\R^{2d})$, $m\in C^\infty_0(\Gamma)$ and assume that the Hessian matrix $\operatorname{Hess}(\phi)$ is non degenerate on $\Gamma$. Then there exists an implicit constant such that for all $\lambda\neq 0$~:
$$ \forall f,g\in \s(\R^d), \qquad \left\|B_\lambda(f,g)\right\|_{L^\infty} \lesssim |\lambda|^{-\delta/2} \|f\|_{L^1} \|g\|_{L^1}.$$
\end{thm}

\dem We repeat the arguments used for Theorem~\ref{thm:thm1curve}. \\
The bilinear kernel of $B_\lambda$ is now given by \be{eq:noyau2} K(x,y,z) = \frac{1}{(2\pi)^d} \int_{\Gamma}
e^{i\eta\cdot (x-y)} e^{i\xi\cdot (x-z)} m(\eta,\xi-\eta) e^{i\lambda\phi(\eta,\xi-\eta)}
d\sigma_\Gamma(\eta,\xi-\eta).\ee Then as $\Gamma$ is a differentiable manifold of dimension $\delta$, up to locally
work, we can use new variables $\eta=(\eta_1,..,\eta_{2d})$ such that $\Gamma$ is described as follows
$$ \Gamma\overset{def}{=}\left\{\eta, \eta_{\delta+1}=...=\eta_{2d}=0 \right\}.$$
The non-degeneracy of the hessian matrix $\operatorname{Hess}(\phi)$ still holds in these new coordinates, so we can apply the stationary phase Lemma on $\Gamma \simeq \R^{\delta}$, which gives
$$ \left|K(x,y,z) \right| \lesssim \lambda^{-\delta/2}.$$
 \findem

\mb Next we want to obtain the appropriate version of Theorem~\ref{thm:thm2smooth}, about $L^2 \times L^1\to L^2$ decays. The proof makes some ``new'' difficulties appear. 

\bigskip

\textit{In the following, we assume that $\Gamma$ can be locally parametrized by any $\delta$-uplet of $(\eta,\xi)\in \R^{2d}$} \footnote{This condition is called the ``non-degeneracy'' of the subspace $\Gamma$ in $\R^{2d}$. This assumption is very important for the study of multilinear operators with symbols admitting singularities on a subspace. We refer the reader to \cite{MTT2} for such results.}. \\

\mb First, let us consider the case of a low dimension $\delta\leq d-1$.
\begin{prop} \label{prop:prop3curve1} In the case where $\delta\leq d-1$, then for all $\lambda\neq 0$, $B_\lambda$ cannot be bounded into $L^2$.
\end{prop}

\dem Indeed, a simple computation gives that $\widehat{B_\lambda(f,g)}$ is a distribution supported on
$$ S\overset{def}{=} \{\eta+\xi, (\eta,\xi)\in \Gamma\}.$$
However as $\Gamma$ is of dimension $\delta$, then $S$ is of dimension less than $\delta$ and so $\textrm{dim}(S)<d$. Then we deduce that
$\widehat{B_\lambda(f,g)}$ is a singular distribution and cannot belongs to $L^2(\R^d)$ in invoking Plancherel Theorem.
\findem

\mb We now deal with the limit case~: $\delta=d$.
\begin{thm} \label{thm:thm3curve1} Assume that $\delta=d$ and that $\Gamma$ can be parametrized by $\eta+\xi$.
Then there exists an implicit constant such that~:
 $$\|B_\lambda (f,g)\|_{L^2} \lesssim \|f\|_{L^1} \|g\|_{L^2}.$$
\end{thm}

\dem Let us just give the sketch of the proof, in applying the same reasoning as for Theorem \ref{thm:thm2smooth}.\\ 
We have to estimate the norm of the operator
\be{int} T_x(f)(\xi)\overset{def}{=} \int_{(\eta,\xi)\in\Gamma} e^{-ix\cdot (\xi-\eta)} e^{i\lambda\phi(\eta,\xi-\eta)} m(\eta,\xi-\eta) f(\eta) d\sigma_{\{\eta, (\eta,\xi-\eta)\in \Gamma\}} (\eta). \ee
We use $\xi\to\eta_{\xi}$ a parametrization such that
$$ (\eta,\xi-\eta)\in\Gamma \Longleftrightarrow \eta=\eta_\xi.$$
The implicit functions theorem permits to deal with such a parametrization and moreover we know that the map $\xi \to \eta_\xi$ is a smooth diffeomorphism. So the integral in (\ref{int}) corresponds to a ``Dirac distribution'' at the point $\eta_\xi$ and we get
$$ T_x(f)(\xi)= e^{i\lambda\phi(\eta_\xi,\xi-\eta_\xi)} m(\eta_\xi,\xi-\eta_\xi) f(\eta_\xi) |\nabla_\xi \eta_\xi|.$$
We also conclude the proof invoking the smooth diffeomorphism $\xi\to \eta_\xi$ and a change of variables.
\findem

\mb Then it remains the more interesting case: $\delta\geq d+1$. \\
If we produce the same reasoning as for Theorem \ref{thm:thm1smooth}, we have to study the linear operator \be{eqTx}
T_x(f)(\xi)\overset{def}{=} \int_{(\eta,\xi)\in\Gamma} e^{-ix\cdot (\xi-\eta)} e^{i\lambda\phi(\eta,\xi-\eta)} m(\eta,\xi-\eta)
f(\eta) d\sigma_{\{\eta,\ (\eta,\xi-\eta)\in\Gamma}(\eta), \ee and then compute the kernel $K$ of $T_xT_x^*$, which gives~:
\begin{align*}
\lefteqn{\int K(\xi,\eta) f(\eta) d\eta =} & & \\
 & &  \int \left(\int f(\eta) \overline{m(\tau,\eta-\tau)} e^{i\lambda(\phi(\tau,\xi-\tau)-\phi(\tau,\eta-\tau))}
d\sigma_{\{\eta,\ (\tau,\eta-\tau)\in\Gamma\}} (\eta)\right) m(\tau,\xi-\tau)
 d\sigma_{\{\tau,\ (\tau,\xi-\tau)\in\Gamma\}}(\tau).
\end{align*}
Then, we would like to compute integrations by parts in the variable $\tau$. The main difficult is that now the integration-domain in the quantity 
$$ \left(\int f(\eta) \overline{m(\tau,\eta-\tau)}d\sigma_{\{\eta,\ (\tau,\eta-\tau)\in\Gamma\}} (\eta)\right) $$
depends on $\tau$ and so will have to be differentiate via the integrations by parts. It is not clear how can we do
this operation, that is why we only consider a particular case.

\gb Assume that $\Gamma$ is an hypersurface as follows~:
$$ \Gamma\overset{def}{=}\left\{ (\eta,\xi)\in\R^{2d},\ \Psi_1(\eta)+\Psi_2(\eta+\xi)=0\right\}$$
where $\Psi_1$ and $\Psi_2$ are smooth functions satisfying
$$ |\nabla \Psi_1|, |\nabla \Psi_2|\geq c$$
for some numerical positive constant $c$.

\begin{thm} \label{thm:thm3curve} Under the above assumptions, with $d\geq 2$, suppose that the derivatives-matrix $\nabla_\xi (\nabla_\eta-\nabla_\xi) \phi$ is non degenerate on
$\operatorname{Supp}(m)$. Then there exists an implicit constant such that for all $\lambda\neq 0$~:
$$ \forall f,g\in \s(\R^d), \qquad \left\|B_\lambda(f,g)\right\|_{L^2} \lesssim |\lambda|^{-d/2} \|f\|_{L^2} \|g\|_{L^1}.$$
\end{thm}

\dem We apply the previous reasoning used for Theorem \ref{thm:thm2smooth}. So we deal with the operator $T_x$ (given
by (\ref{eqTx}))
$$ T_x(f)(\xi)\overset{def}{=} \int e^{-ix\cdot (\xi-\eta)} e^{i\lambda\phi(\eta,\xi-\eta)} m(\eta,\xi-\eta) f(\eta) d\sigma_{\{\eta,\ \Psi_1(\eta)+\Psi_2(\xi)=0\}}(\eta). $$
We compute $T_xT_x^{*}$ and we get~:
\begin{align*}
\lefteqn{T_xT_x^*(f)(\xi)=} & & \\
 & &  \int \left(\int f(\eta) \overline{m(\tau,\eta-\tau)} e^{i\lambda(\phi(\tau,\xi-\tau)-\phi(\tau,\eta-\tau))} d\sigma_{\{\Psi_2(\eta)=\Psi_2(\xi)\}}(\eta)\right) m(\tau,\xi-\tau) d\sigma_{\{\Psi_1(\tau)=-\Psi_2(\xi)\}}(\tau).
\end{align*}
Using integrations by parts in the variable $\tau$, which is feasible as $\{\tau, \Psi_1(\tau)=-\Psi_2(\xi)\}$ is a
manifold around $\xi$ of dimension $d-1\geq 1$, we get~:
$$ \left|T_xT_x^*(f)(\xi)\right| \lesssim \int \int \frac{1}{(1+|\lambda||\xi-\eta|)^N}|f(\eta)|d\sigma_{\{\Psi_2(\eta)=\Psi_2(\xi)\}}(\eta) d\sigma_{\{\Psi_1(\tau)=-\Psi_2(\xi)\}}(\tau)$$
for a large enough integer $N$. Then we divide the space $\R^d$ by unit squares $Q_i\overset{def}{=}i + [0,1]^d$, indexed with
$i\in\Z^d$. So we obtain 
\begin{align*}
\lefteqn{\left\|T_xT_x^*(f)\right\|_{L^2(Q_i)} \lesssim } & & \nonumber \\
 & &  \sum_{j\in\Z^d}
\frac{1}{(1+|\lambda||i-j|)^N} \left\| \int \int_{Q_j} |f(\eta)| d\sigma_{\{\Psi_2(\eta)=\Psi_2(\xi)\}}(\eta)
d\sigma_{\{\Psi_1(\tau)=-\Psi_2(\xi)\}}(\tau) \right\|_{L^2(Q_i)}. 
\end{align*}
For each index $i,j\in\Z^d$, using the
smoothness of the manifolds, we deduce that
\begin{align*}
\lefteqn{\left\| \int \int_{Q_j} |f(\eta)| d\sigma_{\{\Psi_2(\eta)=\Psi_2(\xi)\}}(\eta)
d\sigma_{\{\Psi_1(\tau)=-\Psi_2(\xi)\}}(\tau) \right\|_{L^2(Q_i)}^2} &  \\
& \hspace{1cm} \lesssim \int_{Q_i}\int \int_{Q_j} |f(\eta)|^2
d\sigma_{\{\Psi_2(\eta)=\Psi_2(\xi)\}}(\eta) d\sigma_{\{\Psi_1(\tau)=-\Psi_2(\xi)\}}(\tau) d\xi \\
& \hspace{1cm} \lesssim \int_{0}^\infty \int \int_{Q_j} |f(\eta)|^2 d\sigma_{\{\Psi_2(\eta)=l\}}(\eta)
d\sigma_{\{\Psi_1(\tau)=-l\}}(\tau) \left(\int_{Q_i} \frac{d\sigma_{\Psi_2(\xi)=l}}{|\nabla \Psi_2(\xi)|} \right) dl \\
& \hspace{1cm} \lesssim \int_{Q_j} |f(\eta)|^2 d\eta.
\end{align*}
We have used the ``level-set integration formula'' with respect to $\xi$ and then to $\eta$, the assumptions on the
gradient $\nabla \Psi_2$ and the implicit compactness of the manifolds $\{\eta,\ \Psi_2(\eta)=l\}$
and $\{\tau,\ \Psi_1(\tau)=-\Psi_2(\xi)\}$ due to the compactness support of $m$. By summing in $i,j$, Young's inequality yields
\begin{align}
 \left\|T_xT_x^*(f)\right\|_{2} & \leq \left\| \sum_{i\in\Z^d}\frac{1}{(1+|\lambda||i-j|)^{N}}
\|f\|_{L^2(Q_i)}\right\|_{l^2(j)} \label{eq:ineg} \\
 & \lesssim |\lambda|^{-d} \|f\|_{L^2}. \nonumber
 \end{align}
Remark that (\ref{eq:ineg}) can be seen as a ``off-diagonal estimates''-version of the pointwise bound (\ref{eq:ker}). Then by duality, we conclude that~:
$$ \|T_x\|_{L^2\to L^2} \lesssim |\lambda|^{-d/2}$$
with an implicit constant independent on $x$. \findem

\section{The case of a Coifman-Meyer symbol} \label{sec:3}

We remember that a symbol $m\in L^\infty(\R^{2d})$ is called {\it of Coifman-Meyer type} if it satisfies the bounds
\begin{equation}
\label{CMbound}
\left| \partial_\eta^\alpha \partial_\xi^\beta m(\eta,\xi) \right| \lesssim \frac{1}{\left( |\eta| + |\xi| \right)^{|\alpha| + |\beta|}}
\end{equation}
for sufficiently many multi-indices $\alpha$ and $\beta$.

\gb In the classical work~\cite{coifmanmeyer}, R. Coifman and Y. Meyer show that for such a symbol the bilinear operation $T_{m}$, as defined in~(\ref{defB}), enjoys the same boundedeness properties as given by the H\"older inequality for the standard product (except for extremal values of the Lebesgue indices):
$$\left\| T_{m} (f,g) \right\|_{p} \lesssim \|f\|_q \|g\|_r \quad \textrm{if} \quad \frac{1}{p} = \frac{1}{q} + \frac{1}{r}, \quad 1 < p,q \leq \infty \quad \textrm{and} \quad q <\infty. $$
We now study the boundedness of such bilinear operators by multiplying the symbol $m$ with an extra oscillatory term $e^{i\lambda\phi}$. As previously, we will prove two estimates for
$$B_\lambda(f,g)(x)\overset{def}{=}\frac{1}{(2\pi)^{d/2}}\int e^{ix\cdot (\eta+\xi)} \widehat{f}(\eta) \widehat{g}(\xi) m(\eta,\xi) e^{i\lambda\phi(\eta,\xi)}
d\eta d\xi.$$
The first one describes the decay in $L^\infty$ for $f,g$ belonging to $L^1$ (as for Theorem \ref{thm:thm1smooth}). The second one describes the decay in $L^2$ for $f\in L^2$ and $g\in L^1$ (as for Theorem \ref{thm:thm2smooth}). Then we recover Theorem \ref{thm:thm3smooth} for this kind of symbol by interpolation.

\subsection{The $L^1 \times L^1 \rightarrow L^\infty$ estimate}

\begin{thm} \label{thm:L1}
Let $m$ be a Coifman-Meyer symbol with compact support, and $\phi \in \mathcal{C}^\infty$ such that $\operatorname{Hess} \phi$ is non-degenerate on $\operatorname{Supp} m$. Then for $\lambda\neq 0$
$$
\left\| B_\lambda(f,g) \right\|_\infty \lesssim |\lambda|^{-d} \|f\|_{L^1} \|g\|_{L^1}.
$$
\end{thm}

\dem
The kernel of $B_\lambda$ (by definition: $B_\lambda(f,g)(x) = \int \int K(x,y,z) f(y) g(z)\, dy \,dz$) is given by
$$
K(x,y,z) = \frac{1}{(2\pi)^d} \int e^{i \eta \cdot (x-y)} e^{i \xi \cdot (x-z)} m(\eta,\xi) e^{i\lambda \phi(\eta,\xi)}\, d\eta\, d\xi\,\,.
$$
Of course, it suffices in order to prove the theorem to prove that $K$ is bounded in $L^\infty(\mathbb{R}^{3d})$.

First observe that one can assume $x=y=z=0$: by defining a new phase function $\phi(\eta,\xi) + \frac{1}{\lambda} \eta \cdot (x-y) +\frac{1}{\lambda} \xi \cdot(x-z)$, which we still denote $\phi$, the Hessian remains unchanged. 

We would now like to apply the stationary phase lemma: for the sake of simplicity, assume $\nabla_{\eta,\xi} \phi$ only vanishes at $(\eta_0,\xi_0)$ (the case of several critical points being of course identical), and set $r \overset{def}{=} (\eta_0,\xi_0)$. Next consider a function $\chi$ in $\mathcal{C}^\infty_0$, such that $\chi = 1$ on $B(0,1)$, and $\chi = 0$ on $B(0,2)^c$. Next decompose $K(0,0,0)$ as follows:
\begin{equation}
\begin{split}
K(0,0,0) & = \frac{1}{(2\pi)^d} \int m(\eta,\xi) e^{i\lambda \phi(\eta,\xi)} \left[ \chi \left( \frac{(\xi,\eta)}{\sqrt{\lambda}} \right) + \chi \left( 10 \frac{(\xi,\eta)-(\xi_0-\eta_0)}{r} \right) \right. \\
& \;\;\;\;\;\;\;\;\;\;\;\;\;\;\;\;\;\left. \left[ 1 - \chi \left( \frac{(\xi,\eta)}{\sqrt{\lambda}} \right) - \chi \left( \frac{(\xi,\eta)-(\xi_0-\eta_0)}{10r} \right) \right] \right] d\eta\, d\xi\\
& \overset{def}{=} I + II + III.
\end{split}
\end{equation}
The term $I$ can be estimated brutally:
$$
|I| \lesssim \left| \int \chi \left( \frac{(\xi,\eta)}{\sqrt{\lambda}} \right) d\eta\, d\xi \right| \lesssim \lambda^{-d}.
$$
The stationary phase lemma gives
$$
|II| \lesssim \lambda^{-d}.
$$
Finally, observe that on $\mathbb{R}^{2d} \setminus \left[ B\left( 0,\frac{1}{\sqrt{\lambda}}\right) \cup B \left( (\eta_0,\xi_0), \frac{r}{10} \right) \right]$, $|\nabla \phi| \geq |\eta| + |\xi|$. Thus, integrating by parts $n$ times, with $n=d+2$, gives (we omit the details)
$$
|III| \lesssim  \frac{1}{\lambda^n} \left| \int \left[ 1 - \chi \left( \frac{(\xi,\eta)}{\sqrt{\lambda}} \right) \right] \frac{1}{(|\xi|+|\eta|)^{2n}} d\eta\, d\xi \right| \lesssim \lambda^{-d}.
$$ 
\findem

\mb Using Proposition 5 of Chapter VIII in Stein~\cite{Stein}, we have the associated version of Proposition \ref{prop:prop1smooth}:
\begin{prop}
Let $m$ be a Coifman-Meyer symbol with compact support, and $\phi \in \mathcal{C}^\infty$ such that for a multi-index $\alpha\in\N^{2d}$ $|\alpha|\geq 2$, we have
 $$\left|\partial^\alpha \phi(\eta,\xi) \right|  \gtrsim 1$$ on $\operatorname{Supp}(m)$. Then
$$
\left\| B_\lambda(f,g) \right\|_\infty \lesssim \lambda^{-\frac{1}{2|\alpha|}} \|f\|_{L^1} \|g\|_{L^1}.
$$
\end{prop}

\subsection{The $L^1 \times L^2 \rightarrow L^2$ estimate}

\begin{thm} \label{thm:thm2cm}
Let $m$ be a Coifman-Meyer symbol with compact support, and $\phi \in \mathcal{C}^\infty$ such that $\nabla_\xi \nabla_\eta \phi$ is not singular on $\operatorname{Supp}(m)$. Then
$$
\left\| B(f,g) \right\|_2 \lesssim |\lambda|^{-d/2} \|f\|_{2} \|g\|_{1}.
$$
\end{thm}

\dem First of all, in order to make notations somewhat lighter, we set
$$
\nu(\eta,\xi) \overset{def}{=} m(\eta,\xi-\eta) \;\;\;\;\;\mbox{and}\;\;\;\;\; \Phi(\eta,\xi) = \phi(\eta,\xi-\eta)
$$
and thus for the whole proof of the theorem, $B_\lambda$ will read
$$
B_\lambda (f,g) \overset{def}{=} \mathcal{F}^{-1} \int_{\mathbb{R}^d} \nu(\eta,\xi) e^{i\lambda \Phi(\eta,\xi)} \widehat{f}(\eta) \widehat{g} (\xi-\eta) \,d\eta.
$$
By writing
$$
\mathcal{F} B_\lambda(f,g) (\xi) = \frac{1}{(2\pi)^{d/2}} \int_{\mathbb{R}^d} g(x) \int_{\mathbb{R}^d} e^{-ix\cdot(\xi-\eta)} e^{i\lambda \Phi(\eta,\xi)} \nu(\eta,\xi) \widehat{f}(\eta) \,d\eta \,dx,
$$
and by Plancherel's theorem, it becomes clear that the theorem will follow if one can show that the operator
$$
T_\lambda^x : h \mapsto (Th)(\xi) = \int_{\mathbb{R}^d} e^{-ix\cdot(\xi-\eta)} e^{i\lambda \Phi(\eta,\xi)} \nu(\eta,\xi) h(\eta) \,d\eta,
$$
enjoys the bound
$$
\|T_\lambda^x\|_{L^2 \rightarrow L^2} \lesssim \lambda^{-d/2}.
$$
We now observe that additional hypotheses can be imposed upon $T_\lambda^x$
\begin{itemize}
\item First, due to the non-singularity of $\nabla_\xi \nabla_\eta \Phi$, there holds $|\nabla_\xi \Phi(\eta,\xi) - \nabla_\xi \Phi(\zeta,\xi)| \gtrsim |\zeta-\eta|$ if $\xi \in \operatorname{Supp}\, \nu$ and $|\zeta-\eta| < \epsilon$ for a constant $\epsilon$. By writing $\nu = \nu \sum_j \chi_j$, where the sum is finite and $(\chi_j)$ is a partition of unity such that $\operatorname{Supp} \, \chi_j$ has diameter at most $\frac{1}{2} \epsilon$, one obtains
\begin{equation}
\label{kangourou}
\mbox{if $\zeta, \eta \in \operatorname{Supp}\, \nu$, $|\nabla_\xi \Phi(\eta,\xi) - \nabla_\xi \Phi(\zeta,\xi)| \gtrsim |\zeta-\eta|$}.
\end{equation}
We will assume that this inequality holds in the following.
\item Second, we shall assume that
\begin{equation}
\label{phiphi}
\nabla_\eta \Phi(\eta,0) = 0.
\end{equation}
To see how matters can be reduced to this case, write
$$
(T_\lambda^x h)(\xi) = \int_{\mathbb{R}^d} e^{-ix\cdot(\xi-\eta)} e^{i\lambda \left[ \Phi(\eta,\xi) - \Phi(\eta,0) \right] } \nu(\eta,\xi) \left[ e^{i\lambda \Phi(\eta,0)} h(\eta) \right] \,d\eta,
$$
and notice that $\Phi(\eta,\xi) - \Phi(\eta,0)$ has the desired property, whereas the factor $e^{i\lambda \Phi(\eta,0)}$ multiplying $h(\eta)$ is harmless since bounded.
\item Third, we will suppose that
\begin{equation}
\label{mcondi1}
\nabla_\eta \nu(\eta,0) = 0.
\end{equation}
In order to see why this is possible, consider a Coifman-Meyer symbol $\mu$ such that $\mu = 1$ for $|\xi|\leq |\eta|$. Then write
$$
\nu(\eta,\xi) = \nu(\eta,0) \mu(\eta,\xi) + \left[ \nu(\eta,\xi) - m(\eta,0) \mu(\eta,\xi) \right] \overset{def}{=} \nu_1(\xi,\eta) + \nu_2(\xi,\eta),
$$
and observe that, on the one hand, $\nu_2$ has the desired property; and the other hand, $\nu_1(\eta,\xi)$ is the product of $\mu(\eta,\xi)$, which also has this property, and $\nu(\eta,0)$, which can be directly applied to $h(\eta)$ since $\left\| \nu(\eta,0) h(\eta) \right\|_2 \lesssim \|h(\eta)\|_2$.
\item Finally, we shall assume in the following that
\begin{equation}
\label{mcondi2}
\nu(\eta,\xi)=0 \;\;\;\;\;\;\;\;\;\mbox{if $(\eta,\xi) \in B\left(0,\frac{1}{\sqrt{\lambda}}\right)$}.
\end{equation}
Indeed, select a smooth cut-off function $\chi$ such that $\chi = 1$ on $B\left( 0,\frac{1}{\sqrt{\lambda}}\right)$ and $\chi = 0$ on $B\left( 0,\frac{2}{\sqrt{\lambda}}\right)^c$ and consider the operator $\widetilde{T}^x_\lambda$ with kernel $e^{-ix\cdot(\xi-\eta)} e^{i\lambda \Phi(\eta,\xi)} \nu(\eta,\xi) \chi(\eta,\xi)$. By Sturm's lemma, it enjoys the desired bound on $L^2$:
\begin{equation*}
\begin{split}
\left\|\widetilde{T}^x_\lambda\right\|_{L^2 \rightarrow L^2} & \lesssim \left\| e^{-ix\cdot(\xi-\eta)} e^{i\lambda \Phi(\eta,\xi)} \nu(\eta,\xi) \chi(\eta,\xi) \right\|_{L^\infty_\eta L^1_\xi} \\
& \;\;\;\;\;\;\;\;\;\;\;\;\;\;\;\;\;\;\; + \left\| e^{-ix\cdot(\xi-\eta)} e^{i\lambda \Phi(\eta,\xi)} \nu(\eta,\xi) \chi(\eta,\xi) \right\|_{L^\infty_\xi L^1_\eta} \\
& \lesssim \lambda^{-d/2}.
\end{split}
\end{equation*}
\end{itemize}
Denote $S^x_\lambda$ for the operator
$$
S^x_\lambda = (T^x_\lambda)^* T^x_\lambda
$$
By the classical $T^* T$ argument, $\|T^x_\lambda\|_{L^2 \rightarrow L^2}^2 \leq \|S^x_\lambda\|_{L^2 \rightarrow L^2}$. Thus, in order to prove the theorem it will suffice to show that
$$
\|S^x_\lambda\|_{L^2 \rightarrow L^2} \lesssim \lambda^{-d}
$$
The kernel of $S^x_\lambda$ (by definition $S^x_\lambda h (\eta) = \int K_\lambda^x(\eta,\zeta) h(\zeta) \,d\zeta$) is given by
$$
K_\lambda^x(\eta,\zeta) = e^{ix \cdot(\eta - \zeta)} \int e^{i \lambda(\Phi(\zeta,\xi) - \Phi(\eta,\xi))} \nu(\zeta,\xi) \bar \nu (\eta,\xi) d\xi.
$$
It is clear that the factor $e^{ix \cdot(\eta - \zeta)}$ is irrelevant for the boundedness of $S_\lambda^x$. Thus in the following, we assume $x=0$ and forget about the $x$ superscript.\\
The crucial observation is the following

\begin{prop}
\label{grenouille}
$S_\lambda$ is a singular integral operator. More precisely,
$$
\left\| |\zeta-\eta|^d K_\lambda(\eta,\zeta) \right\|_{L^\infty(\mathbb{R}^{2d})} + \left\| |\eta-\zeta|^{d+1} \nabla_{\eta,\zeta} K_\lambda(\eta,\zeta) \right\|_{L^\infty(\mathbb{R}^{2d})} \lesssim \lambda^{-d}.
$$
\end{prop}
It is natural to try and apply the T1 theorem. We will prove

\begin{prop}
\label{tourterelle}
There holds the bound
$$
\left\| S_\lambda 1 \right\|_{L^\infty} \lesssim  \lambda^{-d}.
$$
\end{prop}

\mb By the T1 theorem of David and Journ\'e~\cite{davidjourne}, the two propositions above, whose proofs follow, give the theorem. \findem

\subsection{Conclusion}

We devote this subsection to deriving some general results from the two previous estimates.

\begin{thm} \label{thm:thm3cm} Let $m$ be a Coifman-Meyer symbol with compact support and $\phi \in \mathcal{C}^\infty$ and assume that $\operatorname{Hess}(\phi)$, $\nabla_\xi (\nabla_\eta-\nabla_\xi) \phi$, $\left(2\nabla_\eta - \nabla_\xi\right)(\nabla_\xi-\nabla_\eta) \phi$ and $\left(\nabla_\eta - 2\nabla_\xi\right)\nabla_\xi \phi$ are non degenerate on $\operatorname{Supp}(m)$. \\
Then for all exponents $p,q,r'\in(1,2]$ satisfying (\ref{eq:condi}), there exists a constant $C=C(p,q,r,\phi,m)$ such that for all $\lambda\neq 0$
$$ \left\|B_\lambda(f,g) \right\|_{L^r} \leq C|\lambda|^{-\frac{d}{2}\left(\frac{1}{p}+\frac{1}{q}-\frac{1}{r}\right)} \|f\|_{L^p} \|g\|_{L^q}.$$
 \end{thm}

\dem We produce the same reasoning as in the case of a smooth symbol (Theorem \ref{thm:thm3smooth} ), one gets a new
set of estimates by interpolating between the $L^1 \times L^1 \rightarrow L^\infty$ and $L^1 \times
L^2 \rightarrow L^2$ cases. \findem

\mb We specify the results when the phase $\phi$ is a homogeneous polynomial function of degree $2$ and a non-compactly
supported symbol $m$:
\begin{thm} \label{thm:thm3cm2} Let $m$ be a Coifman-Meyer symbol and $\phi \in \mathcal{C}^\infty$ a homogeneous polynomial of degree 2, and assume that $\operatorname{Hess}(\phi)$, $\nabla_\xi (\nabla_\eta-\nabla_\xi) \phi$, $\left(2\nabla_\eta - \nabla_\xi\right)(\nabla_\xi-\nabla_\eta) \phi$ and $\left(\nabla_\eta - 2\nabla_\xi\right)\nabla_\xi \phi$ are non degenerate on $\operatorname{Supp}(m)$. \\
Then for all exponents $p,q,r'\in(1,2]$ verifying (\ref{eq:condi}), there exists a constant $C=C(p,q,r,\phi,m)$ such that for all $\lambda\neq 0$
$$ \left\|B_\lambda(f,g) \right\|_{L^r} \leq C|\lambda|^{-\frac{d}{2}\left(\frac{1}{p}+\frac{1}{q}-\frac{1}{r}\right)} \|f\|_{L^p} \|g\|_{L^q}.$$
 \end{thm}

\dem The proof rests on the homogeneity. Let us assume that for $R>>1$, $\gamma_R$ is a smooth  and compactly supported
(on $B(0,2R)$) function such that $\gamma_R=1$ on $B(0,R)\subset \R^{2d}$ and $m_R\overset{def}{=}\gamma_R m$ is still a
Coifman-Meyer symbol. Then it suffices to obtain uniform bound for $B_\lambda^R$ (computed with the truncated symbol
$m_R$) with respect to $R$. We use the scaling as follows: let
$$ \sigma_R(\eta,\xi) = m_R(R\eta,R\xi).$$
So $\sigma_R$ is a uniform Coifman-Meyer symbol and it is supported on $B(0,2)$ and we get by a change of variables~:
$$ B_\lambda^R(f,g)(x)\overset{def}{=} R^{2d} \int e^{iRx\cdot (\eta+\xi)} \widehat{f}(R\eta) \widehat{g}(R\xi) \sigma_R(\eta,\xi) e^{i\lambda\phi(R\eta,R\xi)} d\eta d\xi.$$
Then we use that $\phi$ is a homogeneous polynomial function of order $2$ and so~:
$$ B_\lambda^R(f,g)(x)= R^{2d} \int e^{iRx\cdot (\eta+\xi)} \widehat{f}(R\eta) \widehat{g}(R\xi) \sigma_R(\eta,\xi) e^{i\lambda R^2\phi(\eta,\xi)} d\eta d\xi.$$
We can apply Theorem \ref{thm:thm3cm2} to the symbol $\sigma_R$ and we get~:
\begin{align*}
 \left\|B_\lambda^R(f,g) \right\|_{L^r} & = R^{2d-\frac{d}{r}}\left\| x\to \int e^{ix\cdot (\eta+\xi)} \widehat{f}(R\eta) \widehat{g}(R\xi) \sigma_R(\eta,\xi) e^{i\lambda R^2\phi(\eta,\xi)} d\eta d\xi \right\|_{L^r} \\
& \lesssim R^{2d-\frac{d}{r}} (R^2|\lambda|)^{-\frac{d}{2}\left(\frac{1}{p}+\frac{1}{q}-\frac{1}{r}\right)} R^{-2d} \|f\|_{L^p} R^{\frac{d}{p}} \|g\|_{L^q} R^{\frac{d}{q}} \\
& \lesssim |\lambda|^{-\frac{d}{2}\left(\frac{1}{p}+\frac{1}{q}-\frac{1}{r}\right)} \|f\|_{L^p} \|g\|_{L^q}
\end{align*}
with implicit constants independent on $R$. Then we conclude by passing to the limit when $R\to\infty$.
\findem

\subsection{Proof of Proposition~\ref{grenouille}}

Recall that $K_\lambda$ is given by
$$
K_\lambda(\eta,\zeta) = \int e^{i \lambda(\Phi(\zeta,\xi) - \Phi(\eta,\xi))} \nu(\zeta,\xi) \bar \nu (\eta,\xi) d\xi.
$$

\bigskip

\noindent
\underline{Bound for $\left\| |\zeta-\eta|^d K_\lambda(\zeta,\eta) \right\|_{L^\infty(\mathbb{R}^{2d})}$}

In order to distinguish the cases $|\xi| \geq \frac{1}{\lambda |\zeta - \eta|}$ and $|\xi| \leq \frac{1}{\lambda |\zeta-\eta|}$, we introduce a cutoff function $\chi \in \mathcal{C}^\infty_0$ such that $\chi=1$ on $B(0,1)$ and $\chi = 0$ on $B(0,2)^{c}$, and split the integral defining $K_\lambda$ as follows
\begin{equation}
\label{koala}
\int \dots \,d\xi = \int \chi(\lambda  |\eta-\zeta|\xi) \dots \,d\xi + \int \left[ 1 - \chi(\lambda  |\eta-\zeta|\xi)\right] \dots \,d\xi.
\end{equation}
The first summand in the right hand side of~(\ref{koala}) is estimated directly:
$$
\left| \int \chi(\lambda  |\eta-\zeta|\xi) e^{i \lambda(\Phi(\zeta,\xi) - \Phi(\eta,\xi))} \nu(\zeta,\xi) \bar \nu (\eta,\xi) d\xi \right| \lesssim \left| \int \chi(\lambda |\eta-\zeta| \xi) d\xi \right| \lesssim \lambda^{-d} |\zeta-\eta|^{-d} .
$$
As for the second summand, we will use the identity
\begin{equation}
\label{chat}
M(\xi,\eta,\zeta) \cdot \nabla_\xi e^{i \lambda(\Phi(\zeta,\xi) - \Phi(\eta,\xi))} = e^{i \lambda(\Phi(\zeta,\xi) - \Phi(\eta,\xi))}.
\end{equation}
where
\begin{equation*}
M(\eta,\xi,\zeta) \overset{def}{=} \frac{\lambda \left( \nabla_\xi \Phi(\xi,\eta) - \nabla_\xi \Phi(\xi,\zeta) \right) }{\left| \lambda(\nabla_{\xi}\Phi(\xi,\eta) - \nabla_\xi \Phi(\xi,\zeta) \right|^2}.
\end{equation*}
Notice that, due to~(\ref{kangourou}), $M$ enjoys the bound
\begin{equation}
\label{Mbound}
\left| \nabla^k_{\xi,\eta,\zeta} M^l \right| \lesssim \frac{1}{(\lambda |\zeta-\eta|)^l}\;\;\;\;\;\;\;\mbox{for $k,l \geq 0$}.
\end{equation}
Now integrate by parts in $\xi$ $n$ times (\textit{in the following, $n$ will always denote a big enough integer; for instance, $n = 2d+2$ would suffice everywhere}) using~(\ref{chat}) the second summand of~(\ref{koala}).
When performing these integrations by parts, the derivatives $\nabla_\xi$ might hit either the symbol $\nu (\zeta,\xi) \bar \nu (\eta,\xi)$, or the cut-off function $\left[ 1 - \chi(\lambda  |\eta-\zeta|\xi)\right]$, or $M(\eta,\xi)$. To simplify the notations, we simply consider the situation where all the derivatives hit one of these three factors; furthermore, we will a bit abuse notations by not keeping track of the vectorial relations, rather treating all the factors as scalars.
\begin{itemize}
\item If it is $M(\eta,\xi)$, use the bound~(\ref{Mbound}) to obtain
\begin{equation*}
\begin{split}
&\left| \int e^{i \lambda(\Phi(\zeta,\xi) - \Phi(\eta,\xi))}  \nabla^n_\xi M(\eta,\xi)^n \left[ \nu(\zeta,\xi) \bar \nu (\eta,\xi) \right] \left[ 1 - \chi(\lambda |\eta-\zeta| \xi)\right] d\xi \right| \\
&\;\;\;\;\;\;\;\;\;\;\;\;\;\;\;\;\;\;\;\;\;\;\;\;\;\;\; \lesssim \frac{1}{(\lambda|\zeta-\eta|)^n}.
\end{split}
\end{equation*}
\item If it is $\nu(\zeta,\xi) \bar \nu (\eta,\xi)$, use furthermore that the symbol satisfies~(\ref{CMbound}) to get
\begin{equation*}
\begin{split}
& \left| \int e^{i \lambda(\Phi(\zeta,\xi) - \Phi(\eta,\xi))} M(\eta,\xi)^n \nabla^n_\xi \left[ \nu(\zeta,\xi) \bar \nu (\eta,\xi) \right] \left[ 1 - \chi(\lambda |\eta-\zeta| \xi)\right] d\xi \right| \\
&\;\;\;\;\;\;\;\;\;\;\;\;\;\;\;\;\;\;\;\;\;\;\;\;\;\;\; \lesssim \frac{1}{(\lambda|\zeta-\eta|)^n} \int_{r \gtrsim \frac{1}{\lambda |\zeta - \eta|}} \frac{1}{r^n} r^{d-1} dr \lesssim \lambda^{-d} |\zeta-\eta|^{-d}.
\end{split}
\end{equation*}
\item Finally, if it is $\left[ 1 - \chi(\lambda |\eta-\zeta| \xi)\right]$, one gets
\begin{equation*}
\begin{split}
& \left| \int e^{i \lambda(\Phi(\zeta,\xi) - \Phi(\eta,\xi))} M(\eta,\xi)^n \nabla^n_\xi \left[ 1 - \chi(\lambda |\eta-\zeta| \xi)\right] \nu(\zeta,\xi) \bar \nu (\eta,\xi) d\xi \right| \\
&\;\;\;\;\;\;\;\;\;\;\;\;\;\;\;\;\;\;\;\;\;\;\;\;\;\;\; \lesssim \frac{1}{(\lambda|\zeta-\eta|)^n} \int_{r \sim \frac{1}{\lambda |\zeta - \eta|}} \lambda^n |\zeta-\eta|^n r^{d-1} dr \lesssim \lambda^{-d} |\zeta-\eta|^{-d}.
\end{split}
\end{equation*}
\end{itemize}
The above estimates give
$$
K_\lambda(\eta,\zeta) \lesssim \lambda^{-d} |\zeta-\eta|^{-d} + \lambda^{-n} |\zeta-\eta|^{-n}.
$$
Since obviously $|K_\lambda (\eta,\zeta)| \lesssim 1$, the proposition follows. \findem

\bigskip

\noindent
\underline{Bound for $\left\| |\zeta-\eta|^{d+1} \nabla K_\lambda(\zeta,\eta) \right\|_{L^\infty(\mathbb{R}^{2d})}$}

Due to the symmetry of $K$, it suffices to bound $\nabla_\eta K$. Applying $\nabla_\eta$ to $K_\lambda(\zeta,\eta)$ yields
\begin{equation*}
\begin{split}
\nabla_\eta K_\lambda(\zeta,\eta) & = \int i \lambda \nabla_\eta \Phi(\eta,\xi) e^{i \lambda(\Phi(\zeta,\xi) - \Phi(\eta,\xi))} \nu(\zeta,\xi) \bar \nu (\eta,\xi) d\xi \\
& \;\;\;\;\;\;\;\;\;\;\;\;\;\;\;\;\;\;\;\;\;\;\; + \int  e^{i \lambda(\Phi(\zeta,\xi) - \Phi(\eta,\xi))}  \nu(\zeta,\xi)\nabla_\eta \bar \nu (\eta,\xi) d\xi \\
\end{split}
\end{equation*}

The key observation is that, denoting $\mu(\xi,\eta,\zeta)$ for either
$$
\lambda \nabla_\eta \Phi(\eta,\xi) \nu(\zeta,\xi) \bar \nu (\eta,\xi)\;\;\;\;\;\mbox{or}\;\;\;\;\; \nu(\zeta,\xi)\nabla_\eta \bar \nu (\eta,\xi),
$$
its $\xi$ derivatives enjoy the pointwise bound
\begin{equation}
\label{mubound}
\mbox{on $\textit{Supp} \nu$,} \;\;\;\;\;\left| \nabla_\xi^{k} \mu(\xi,\eta,\zeta) \right| \lambda |\xi|^{1-k}.
\end{equation}

This follows from the following facts:
\begin{itemize}
\item On the one hand, the assumption~(\ref{phiphi}) implies $|\nabla_\eta \Phi(\eta,\xi) | \lesssim |\xi|$.
\item On the other hand, the assumptions~(\ref{mcondi1}),~(\ref{mcondi2}) and~(\ref{CMbound}) give $\left| \nabla^k_\xi \nabla_\eta \nu(\eta,\xi) \right| \lesssim \lambda |\xi|^{1-k}$.
\end{itemize}
Thus it suffices to prove that, for $\mu$ satisfying the above bound,
$$
\left| \int \mu(\xi,\eta,\zeta) e^{i \lambda(\Phi(\zeta,\xi) - \Phi(\eta,\xi))}\,d\xi \right| \lesssim \lambda^{-d/2}.
$$
Using a smooth cut-off function, we now split the above integral into two, with integration domains respectively the regions $|\xi| \geq \frac{1}{\lambda |\zeta - \eta|}$ and $|\xi| \leq \frac{1}{\lambda |\zeta-\eta|}$ \textit{However, for the remainder of this article, and for the sake of simplicity in the notations, cut-off functions will not appear explicitly and we will simply write}
\begin{equation}
\label{eqI}
\int \mu(\xi,\eta,\zeta) e^{i \lambda(\Phi(\zeta,\xi) - \Phi(\eta,\xi))}\,d\xi = \int \dots d\xi = \int_{|\xi|\lesssim \frac{1}{\lambda|\zeta-\eta|}} \dots d\xi + \int_{|\xi| \gtrsim \frac{1}{\lambda|\zeta-\eta|}} \dots d\xi.
\end{equation}
\textit{Similarly, we will not care in the estimates when derivatives hit the cutoff function: it should be clear from above that this always produces harmless terms.} \\
The first summand in~(\ref{eqI}) can be dealt with exactly as in the estimate of $\left\| |\zeta-\eta|^d K_\lambda(\zeta,\eta) \right\|_{L^\infty(\mathbb{R}^{2d})}$.

\mb For the second summand in~(\ref{eqI}), also proceed as in the estimate of $\left\| |\zeta-\eta|^d K_\lambda(\zeta,\eta) \right\|_{L^\infty(\mathbb{R}^{2d})}$, namely integrate by parts $n$ times using the identity~(\ref{chat}). Just like there, the worst term here occurs when $\nabla_\xi$ hits the ``symbol'' $\mu(\zeta,\xi,\eta)$, namely
$$
\int_{|\xi|\gtrsim \frac{1}{\lambda|\zeta-\eta|}} M(\xi,\eta)^n e^{i \lambda(\Phi(\zeta,\xi) - \Phi(\eta,\xi))} \nabla^n_\xi \mu(\xi,\eta,\zeta) \,d\xi.
$$
But, due to the bound~(\ref{mubound}), it is easily estimated by
$$
\frac{C}{(\lambda |\zeta-\eta|)^n} \int_{r \gtrsim \frac{1}{\lambda |\zeta-\eta|}} \lambda r^{d-n}\,dr \lesssim \lambda^{-d} |\zeta-\eta|^{-d-1},
$$
which is the desired result. \findem

\subsection{Proof of proposition~\ref{tourterelle}}

We want to prove that
$$
\left[ S_\lambda 1 \right] (\eta) = \int \int e^{i \lambda(\Phi(\zeta,\xi) - \Phi(\eta,\xi))} \nu(\zeta,\xi) \bar \nu (\eta,\xi) \,d\xi \,d\zeta.
$$
belongs to $L^\infty$, with the bound
$$
\left\| S_\lambda 1 \right\|_{\infty} \lesssim \lambda^{-d}.
$$
This will be achieved by splitting the integral into several pieces corresponding to different integration domains, and
estimating them separately. As above, this is done with the help of cut-off functions, and we adopt the same convention
that they will not appear explicitly.

Of course, the whole idea is to take advantage of the oscillations by integrating by parts
\begin{itemize}
\item Either in $\xi$, and we rely then on~(\ref{chat}) and~(\ref{Mbound}).
\item Or in $\zeta$; but what are the critical points of $\Phi(\xi,\zeta)$ in $\zeta$? We know that any $(\zeta,0)$ is one, by the assumption~(\ref{phiphi}). There may be other ones, but due to the hypothesis that $\nabla_\xi \nabla_\zeta \Phi$ is invertible, they can occur only away from the plane $\{ \xi = 0 \}$. Suppose that for some $\xi_0,\zeta_0$, $\nabla_\zeta \Phi (\xi_0,\zeta_0)=0$. By the implicit function theorem, and invertibility of $\nabla_\xi \nabla_\zeta \Phi$, there is a smooth $d$-dimensional surface on which $\nabla_\zeta \Phi$ vanishes. We shall however consider in the following that $\nabla_\zeta \Phi$ vanishes only on $\{ \xi = 0 \}$. This is simply because the possible other singularity planes are easier to treat, since the function $\nu(\zeta,\xi) \bar \nu (\eta,\xi)$ is most singular for $\xi=0$. Thus we shall assume that
\begin{equation}
\label{singe}
\mbox{on $\operatorname{Supp}\, \nu$,}\;\;\; |\nabla_\zeta \Phi (\zeta,\xi)| \gtrsim |\xi|.
\end{equation}
This will be used via the identity
\begin{equation}
\label{herisson}
N(\xi,\eta,\zeta) \cdot \nabla_\zeta e^{i \lambda(\Phi(\zeta,\xi) - \Phi(\eta,\xi))} = e^{i \lambda(\Phi(\zeta,\xi) - \Phi(\eta,\xi))}.
\end{equation}
where
\begin{equation*}
N(\xi,\eta,\zeta) \overset{def}{=} \frac{\lambda  \nabla_\zeta \Phi(\xi,\zeta) }{\left| \lambda(\nabla_{\zeta}\Phi(\xi,\zeta)\right|^2}.
\end{equation*}
enjoys the bound
\begin{equation}
\label{Nbound}
\left| \nabla^k_{\xi,\eta,\zeta} N^l(\xi,\eta,\zeta) \right| \lesssim \frac{1}{(\lambda |\xi|)^l} \frac{1}{|\xi|^k} \;\;\;\;\;\;\;\mbox{for $k,l \geq 0$}.
\end{equation}
\end{itemize}

\bigskip
\noindent
\underline{$L^\infty$ bound for the piece $|\zeta|\sim |\eta|$ and $|\xi| \lesssim |\eta|$.}
We are considering
$$
F_1 (\eta) \overset{def}{=} \int \int_{\substack{|\zeta|\sim |\eta|\\|\xi| \lesssim |\eta|}}
e^{i \lambda(\Phi(\zeta,\xi) - \Phi(\eta,\xi))} \nu(\zeta,\xi) \bar \nu (\eta,\xi) \,d\xi \,d\zeta
$$
Fix $\eta_0$. By~(\ref{kangourou}) and~(\ref{singe}), the phase $\Phi(\zeta,\xi) - \Phi(\eta_0,\xi)$ is stationary for $\zeta=\eta_0$, $\xi=0$. The Hessian is given by
$$
\operatorname{Hess}_{\xi,\zeta}\left[ \Phi(\zeta,\xi) - \Phi(\xi,\eta_0) \right] (0,\eta_0) = \left( \begin{array}{cc} 0 & \nabla_\xi \nabla_\zeta \Phi(\xi,\zeta) \\ \nabla_\xi \nabla_\eta \Phi(\xi,\zeta) & \nabla_\zeta \nabla_\zeta \Phi(\xi,\zeta) \end{array} \right),
$$
hence it is invertible by hypothesis. An application of the stationary phase principle gives thus
$$
\left| F_1 (\eta_0) \right| \lesssim \lambda^{-d}.
$$

\bigskip
\noindent
\underline{$L^\infty$ bound for the piece $|\zeta|\sim |\eta|$, $|\xi| \gtrsim |\eta|$ if $|\eta| \gtrsim \frac{1}{\sqrt{\lambda}}$.}
We are now considering
$$
F_2 (\eta) \overset{def}{=} \int \int_{\substack{|\zeta|\sim |\eta|\\|\xi| \gtrsim |\eta|}}
e^{i \lambda(\Phi(\zeta,\xi) - \Phi(\eta,\xi))} \nu(\zeta,\xi) \bar \nu (\eta,\xi) \,d\xi \,d\zeta.
$$
Integrate by parts $n$ times in $\zeta$ using~(\ref{herisson}). The worst terms occur when all the $\zeta$ derivatives hit $\nu(\zeta,\xi)$ or $N(\xi,\eta,\zeta)^n$:
$$
\int \int_{\substack{|\zeta|\sim |\eta|\\|\xi| \gtrsim |\eta|}} \nabla^{n-\ell}_\zeta N(\xi,\eta,\zeta)^n
e^{i \lambda(\Phi(\zeta,\xi) - \Phi(\eta,\xi))} \nabla^\ell_\zeta \nu(\zeta,\xi) \bar \nu (\eta,\xi) \,d\xi \,d\zeta \;\;\;\;\mbox{for $0\leq \ell \leq n$}.
$$
By~(\ref{CMbound}) and~(\ref{Nbound}), this can be bounded by
$$
\int \int_{\substack{|\zeta|\sim |\eta|\\|\xi| \gtrsim |\eta|}} \frac{1}{(\lambda |\xi|)^n} \frac{1}{|\xi|^n}\,d\xi \,d\zeta \lesssim \lambda^{-n} |\eta|^{d-2n} \lesssim \lambda^{-d},
$$
where the last inequality follows from the assumption that $|\eta| \gtrsim \frac{1}{\sqrt{\lambda}}$. This is the desired estimate.

\bigskip
\noindent
\underline{$L^\infty$ bound for the piece $|\zeta| << |\eta|$, $|\zeta| + |\xi| \gtrsim \frac{1}{\lambda |\eta|}$ if $|\eta| \gtrsim \frac{1}{\sqrt{\lambda}}$.}
This corresponds to
$$
F_3 (\eta) \overset{def}{=} \int \int_{\substack{|\zeta| << |\eta|\\|\zeta| + |\xi| \gtrsim \frac{1}{\lambda |\eta|}}}
e^{i \lambda(\Phi(\zeta,\xi) - \Phi(\eta,\xi))} \nu(\zeta,\xi) \bar \nu (\eta,\xi) \,d\xi \,d\zeta.
$$
Integrate by parts $n$ times in $\xi$ with the help of~(\ref{chat}). As usual, the worst term is
$$
\int \int_{\substack{|\zeta| << |\eta|\\|\zeta| + |\xi| \gtrsim \frac{1}{\lambda |\eta|}}} M(\xi,\eta,\zeta)^n
e^{i \lambda(\Phi(\zeta,\xi) - \Phi(\eta,\xi))} \nabla_\xi^n \left[ \nu(\zeta,\xi) \bar \nu (\eta,\xi) \right] \,d\xi \,d\zeta.
$$
By~(\ref{CMbound}), (\ref{Mbound}), and the restriction on the integration domain which implies $|\zeta-\eta| \gtrsim |\eta|$, this can be bounded by
$$
\left| \frac{1}{(\lambda |\eta|)^n} \int \int_{\substack{|\zeta| << |\eta|\\ |\zeta| + |\xi| \gtrsim \frac{1}{\lambda |\eta|}}} \left( \frac{1}{(|\zeta|+|\xi|)^n} + \frac{1}{(|\eta|+|\xi|)^n} \right) \,d\xi \,d\zeta \right|
\lesssim \lambda^{-d},
$$
where the last inequality follows from the assumption that $|\eta| \gtrsim \frac{1}{\sqrt{\lambda}}$.

\bigskip
\noindent
\underline{$L^\infty$ bound for the piece $|\zeta|<< |\eta|$, $|\zeta| + |\xi| \lesssim << \frac{1}{\lambda |\eta|}$, if $|\eta| \gtrsim \frac{1}{\lambda}$}
This corresponds to
$$
F_4 (\eta) \overset{def}{=} \int \int_{\substack{|\zeta|<< |\eta|\\|\zeta| + |\xi| \lesssim \frac{1}{\lambda |\eta| }} }
e^{i \lambda(\Phi(\zeta,\xi) - \Phi(\eta,\xi))} \nu(\zeta,\xi) \bar \nu (\eta,\xi) \,d\xi \,d\zeta.
$$
A direct estimate gives
$$
\left| F_4 (\eta) \right| \lesssim (\lambda |\eta|)^{-2d} \lesssim \lambda^{-d},
$$
where the last inequality is justified by the hypothesis that $|\eta| \gtrsim \frac{1}{\lambda}$.

\bigskip
\noindent
\underline{$L^\infty$ bound for the piece $|\zeta|>>|\eta|$ if $|\eta| \gtrsim \frac{1}{\sqrt{\lambda}}$}
We are now considering
$$
F_5 (\eta) \overset{def}{=} \int \int_{|\zeta| >> |\eta|}
e^{i \lambda(\Phi(\zeta,\xi) - \Phi(\eta,\xi))} \nu(\zeta,\xi) \bar \nu (\eta,\xi) \,d\xi \,d\zeta.
$$
Integrate by parts $n$ times in $\xi$; the worst resulting term is
$$
\int \int_{|\zeta| >> |\eta|} M(\xi,\eta,\zeta)^n
e^{i \lambda(\Phi(\zeta,\xi) - \Phi(\eta,\xi))} \nabla^n_\xi \left[ \nu(\zeta,\xi) \bar \nu (\eta,\xi) \right] \,d\xi \,d\zeta.
$$
By~(\ref{CMbound}),~(\ref{Mbound}), and since the restrictions on the integration domain imply $|\eta - \zeta| \gtrsim |\zeta|$, this can be bounded by
$$
\int_{|\zeta| >> |\eta|} \frac{1}{(\lambda|\zeta|)^n} \left[ \frac{1}{(|\zeta|+|\xi|)^n} + \frac{1}{(|\eta|+|\xi|)^n} \right] \,d\xi\,d\zeta \lesssim \lambda^{-n}|\eta|^{2d-2n} \lesssim \lambda^{-d},
$$
where the last inequality follows from the assumption $|\eta|> \frac{1}{\sqrt{\lambda}}$.

\bigskip
\noindent
\underline{$L^\infty$ bound for the piece $|\zeta| \gtrsim \frac{1}{\sqrt{\lambda}}$ and $|\xi| \gtrsim \frac{1}{\sqrt{\lambda}}$ if $|\eta| << \frac{1}{\sqrt{\lambda}}$}
This corresponds to
$$
F_6 (\eta) \overset{def}{=} \int \int_{\substack{|\zeta| \gtrsim \frac{1}{\sqrt{\lambda}} \\ |\xi| \gtrsim \frac{1}{\sqrt{\lambda}}}}
e^{i \lambda(\Phi(\zeta,\xi) - \Phi(\eta,\xi))} \nu(\zeta,\xi) \bar \nu (\eta,\xi) \,d\xi \,d\zeta.
$$
After $n$ integrations by parts in $\xi$ (using~(\ref{chat}), and $n$ integrations by parts in $\zeta$ (using~(\ref{herisson}), the worst terms are
$$
\int \int_{\substack{|\zeta| \gtrsim \frac{1}{\sqrt{\lambda}} \\ |\xi| \gtrsim \frac{1}{\sqrt{\lambda}}}}
M(\xi,\eta,\zeta) \nabla^k_\xi \nabla_\zeta^\ell N(\xi,\eta,\zeta) e^{i \lambda(\Phi(\zeta,\xi) - \Phi(\eta,\xi))}
\nabla^{n-k}_\xi \nabla^{n-\ell}_\zeta \left[ \nu(\zeta,\xi) \bar \nu (\eta,\xi) \right] \,d\xi \,d\zeta
$$
for $0\leq k,\ell \leq n$.
This can be bounded by
$$
\int \int_{\substack{|\zeta| \gtrsim \frac{1}{\sqrt{\lambda}} \\ |\xi| \gtrsim \frac{1}{\sqrt{\lambda}}}}
\frac{1}{(\lambda|\zeta|)^n} \frac{1}{(\lambda|\xi|)^n} \frac{1}{|\xi|^{k+\ell}} \frac{1}{|\xi|^{2n-k-\ell}}.
$$
\bigskip
\noindent
\underline{$L^\infty$ bound for the piece $|\zeta| \lesssim \frac{1}{\sqrt{\lambda}}$ and $|\xi| \gtrsim \frac{1}{\sqrt{\lambda}}$ if $|\eta| << \frac{1}{\sqrt{\lambda}}$}
The function under consideration is now
$$
F_7 (\eta) \overset{def}{=} \int \int_{\substack{|\zeta| \lesssim \frac{1}{\sqrt{\lambda}} \\ |\xi| \gtrsim \frac{1}{\sqrt{\lambda}}}}
e^{i \lambda(\Phi(\zeta,\xi) - \Phi(\eta,\xi))} \nu(\zeta,\xi) \bar \nu (\eta,\xi) \,d\xi \,d\zeta.
$$
It can be estimated by integrating by parts $n$ times in $\zeta$; the details are left to the reader.

\bigskip
\noindent
\underline{$L^\infty$ bound for the piece $|\xi|<<\frac{1}{\sqrt{\lambda}}$ if $|\eta|<<\frac{1}{\sqrt{\lambda}}$} In this case, $\nu
(\xi,\eta)=0$ by~(\ref{mcondi2}), thus the function we want to bound is identically zero!

\bigskip

The previous estimates yield the desired conclusion $\|S_\lambda 1\|_\infty \lesssim \lambda^{-d}$. Indeed, if $|\eta| \gtrsim \frac{1}{\sqrt{\lambda}}$, $S_\lambda 1 = F_1 + \dots + F_5$; and if $|\eta| << \frac{1}{\sqrt{\lambda}}$, $S_\lambda 1 = F_6 + F_7$. \findem

\section{The case of an $x$-dependent symbol} \label{sec:4}

We devote this section to results concerning our bilinear oscillatory integrals, involving an $x$-dependent symbol $m$. Firstly we extend the previous results about Coifman-Meyer symbols.

\begin{prop} Assume that $\phi$ satisfies the above assumptions and let $m$ be an $x$-dependent Coifman-Meyer symbol~:
\be{eq:nonh} \left| \partial_\eta^\alpha \partial_\xi^\beta \partial_x^\gamma  m(x,\eta,\xi) \right| \lesssim \frac{1}{\left( 1+  |\eta| + |\xi| \right)^{|\alpha| + |\beta|}} \ee
 for sufficiently many indices $\alpha,\beta$ and $\gamma$. \\
Then for all exponents $p,q,r\in(1,2]$ verifying (\ref{eq:condi}) and any power $exp$ satisfying 
$$ exp\geq -\frac{d}{2} \left(\frac{1}{p}+\frac{1}{q}-\frac{1}{r}\right) \quad \textrm{ and } \quad exp>-d+\frac{1}{p}+\frac{1}{q}$$
there exists a constant $C=C(p,q,r,\phi,m)$ such that for all $|\lambda|\geq 1$
$$ \left\|B_\lambda(f,g) \right\|_{L^r} \leq C |\lambda|^{exp}\|f\|_{L^p} \|g\|_{L^q}.$$
\end{prop}

\mb For $x$-dependent symbols, we cannot keep the same decay (than for $x$-independent symbols) relatively to $|\lambda|$ due to some extra integrations by parts, however we get some boundedness in the product of Lebesgue spaces for these new bilinear operators $B_1$.

\mb We emphasize that the proof of the $L^1\times L^1 \to L^\infty$ decay (see Theorem \ref{thm:L1}) still holds for an $x$-dependent symbol $m$. So this proposition is only interesting for the other exponents.

\mb
\dem We will follow ideas of Theorem 34 in \cite{coifmanmeyer}, which were extended in a more general framework by the first author in Section 4 of \cite{B2}. We just explain the main ideas and leave the details to the reader. This reasoning permits to reduce the study of $x$-dependent bilinear symbols to the one of $x$-independent symbols. The main tool is some ``off-diagonal decay'', which is stronger than the global boundedness estimate. \\
In fact, we claim that for an $x$-independent symbols $m$ verifying the non-homogeneous decay (\ref{eq:nonh}) and for all square $I\subset \R^d$ of measure $1$~:
\be{eq:off} \|B_\lambda(f,g)\|_{r,I} \lesssim |\lambda|^{e}
\left[\sum_{k\geq 0} 2^{-k(\frac{1}{p}+\epsilon)} \|f\|_{p,2^{k+1}I} \right] \left[ \sum_{k\geq 0} 2^{-k(\frac{1}{q}+\epsilon)} \|g\|_{q,2^{k+1}I}
\right],
\ee
for every $\epsilon>0$ with 
$$e=\max\left\{-\frac{d}{2} \left(\frac{1}{p}+\frac{1}{q}-\frac{1}{r}\right),\ -d+\frac{1}{p}+\frac{1}{q}+2\epsilon\right\}.$$
For $J$ an interval, we denote by $\|\ \|_{p,J}$ the $L^p$ norm on $J$ and for a positive real $\lambda$ by $\lambda J$ the interval of lenght $\lambda |J|$ with the same center than $J$. \\ 
(\ref{eq:off}) comes from the fast decay of the bilinear kernel $K$ of $B_\lambda$ away from the diagonal~:
$$ \left| K(x,y,z)\right| = \left| \int_{\R^{2d}} e^{-i\left[(x-y)\cdot \eta + (x-z)\cdot \xi \right]} e^{i\lambda \phi(\eta,\xi)} m(\eta,\xi)\,d\eta \,d\xi \right| \lesssim  |\lambda|^{-d}\left(1+\frac{|x-y|}{|\lambda|}+\frac{|x-z|}{|\lambda|}\right)^{-N}$$
for all integer $N\geq 0$ by computing integrations by parts and using (\ref{eq:nonh}).
So in decomposing $f$ and $g$ on the dyadic coronas around $I$, the terms for $k>1$ in (\ref{eq:off}) come from easily as $p,q\geq 1$ 
with the exponent $e=-d+\frac{1}{p}+\frac{1}{q}+2\epsilon$ and the term for $k=0$ is due to the previous theorem. \\
Then from (\ref{eq:off}) for $x$-independent symbols, we deduce (\ref{eq:off}) for $x$-dependent symbols $m$ in using a Sobolev imbedding and in considering the space variable $x$ of $m$ independently to the variable $x$ of $B_\lambda(f,g)$ (see Lemma 6, Chap. VI of \cite{coifmanmeyer} and Theorem 4.5 of \cite{B2} for similar arguments).
So let us take the $x$-dependent symbol $m$ of the statement. Since (\ref{eq:off}) holds for each square $I$ of measure $1$, we use a bounded covering of $\R^{d}$ by such squares $(I_i)_i$. Denoting $r_{p,q}$ the exponent satisfying~:
$$ \frac{1}{r_{p,q}} = \frac{1}{p}+\frac{1}{q},$$
we remark that $r\geq r_{p,q}$ from (\ref{eq:condi}). So
$$ \|B_\lambda (f,g)\|_{r} \simeq \left\| \|B_\lambda(f,g)\|_{r,I_i} \right\|_{l^r(i)} \lesssim \left\| \|B_\lambda(f,g)\|_{r,I_i} \right\|_{l^{r_{p,q}}(i)}.$$
We use (\ref{eq:off}) for all $i$ and then H\"older inequality permits to deduce the desired inequality. \findem

\mb Then we would like to describe a more general estimate than the one describing the $L^1 \times L^1 \to L^\infty$ boundedness in Theorem \ref{thm:L1}. More precisely, we know that a $x$-independent Coifman-Meyer symbol $m$ yields a $2d$-dimensionnal linear (convolution) Calder\'on-Zygmund kernel $K(y,z)$ by
$$ K(y,z):= \int_{\R^{2d}} e^{i(\eta\cdot y + \xi\cdot z)} m(\eta,\xi) d\eta d\xi.$$
By the well-known theory of linear Calder\'on-Zygmund operators, the corresponding linear operator is bounded from the Hardy space $\mathcal{H}^1(\R^{2d})$ into $L^1(\R^{2d})$.
So we are now interested by a symbol $m$ obtained as the Fourier transform of a Calder\'on-Zygmund kernel. This case is also more general than the Coifman-Meyer case. We obtain a weaker version of the $L^1\times L^1\to L^\infty$ estimate in this case using the Hardy space.

\begin{thm}
Let $m$ be a bounded (non smooth) compactly supported symbol such that the distributional kernel $K_m$ defined by
$$ K_m(x,y,z) \overset{def}{=} \int_{\R^{2d}} e^{i(\eta\cdot y + \xi\cdot y)} m(x,\eta,\xi) d\eta d\xi$$
satisfies for all fixed $x$ 
$$ \left|\partial_y^a \partial_z^b K_m(x,y,z)\right|\lesssim \frac{1}{(|y|+|z|)^{2d+a+b}}$$
for every multi-index $a,b\in \N^{2d}$ with $|a|,|b|\leq d+1$.
Let $\phi \in \mathcal{C}^\infty$ such that $\operatorname{Hess} \phi$ is non-degenerate on $\operatorname{Supp} m$. Then for $\lambda\neq 0$
$$
\left\| B_\lambda(f,g) \right\|_\infty \lesssim |\lambda|^{-d} \|f\|_{\mathcal{H}^1} \|g\|_{\mathcal{H}^1}.
$$
\end{thm}

\dem \underline{Step 1: bound on the kernel in BMO}

The kernel of $B_\lambda$ (by definition: $B_\lambda(f,g)(x) = \int \int K(x,y,z) f(y) g(z)\, dy \,dz$) is given by
$$
K(x,y,z) \overset{def}{=} \frac{1}{(2\pi)^d} \int e^{i \eta \cdot (x-y)} e^{i \xi \cdot (x-z)} m(x,\eta,\xi) e^{i\lambda \phi(\eta,\xi)} d\eta d\xi\,\,.
$$
Take a function $\chi \in \mathcal{C}_0^\infty$ such that $\chi = 1$ on $\operatorname{Supp} m$.
Let us fix the point $x=x_0$ and denote by $M_{x_0}$ the linear operator on $\mathbb{R}^{2d}$ with symbol $m(x_0,\cdot)$ ($M_{x_0}$ corresponds to the convolution in $\R^{2d}$ by $(y,z)\to K_m(x_0,y,z)$), and by $F$ the function whose Fourier transform (in $\mathbb{R}^{2d}$) is $\widehat{F}(\eta,\xi) = \chi(\eta,\xi) e^{i\lambda \phi}(\eta,\xi)$. It is then possible to write
$$
K(x_0,y,z) = (M_{x_0}F)(x_0-y,x_0-z).
$$
On the one hand, by stationary phase, $\|F\|_\infty \lesssim \lambda^{-d}$; on the other hand, by standard properties of Calder\'on-Zygmund operators, the operator $M_{x_0}$ is bounded from $L^\infty(\mathbb{R}^{2d})$ to $BMO(\mathbb{R}^{2d})$ (as the assumptions imply that $M$ is a Calder\'on-Zygmund operator on $\R^{2d}$). Therefore,
\begin{equation}
\label{ecureuil}
\left\| M_{x_0}F \right\|_{BMO} \lesssim \lambda^{-d}.
\end{equation}

\bigskip

\noindent \underline{Step 2: duality between $\mathcal{H}^1$ and $BMO$ argument}

We first need to rewrite $B_\lambda(f,g)$: if $x_0 \in \mathbb{R}^d$,
\begin{equation}
\begin{split}
B_\lambda(f,g)(x_0) & = \int K(x_0,y,z) f(y) g(z) \,dy\,dz \\
& = \int (M_{x_0}F)(u,v) f(x_0-u) g(x_0-v) \,du \,dv \\
& = \langle (M_{x_0}F)(u,v) \,,\,f(x_0-u) g(x_0-v) \rangle\,\,.
\end{split}
\end{equation}
Using first duality between $\mathcal{H}^1(\R^{2d})$ and $BMO(\R^{2d})$, and then the estimate~(\ref{ecureuil}) gives
$$
|B_\lambda(f,g)(x_0)| \leq \|M_{x_0}F\|_{BMO(\mathbb{R}^{2d})} \left\| f(x_0-u) g(x_0-v) \right\|_{\mathcal{H}^1(\mathbb{R}^{2d})} \lesssim \lambda^{-d} \left\| f(x_0-u) g(x_0-v) \right\|_{\mathcal{H}^1(\mathbb{R}^{2d})}
$$
Using finally Lemma \ref{lemme} (and invariance of Hardy spaces under translations), we get the desired result, namely
$$
|B_\lambda(f,g)(x_0)| \lesssim \lambda^{-d} \|f\|_{\mathcal{H}^1(\mathbb{R}^{d})} \|g\|_{\mathcal{H}^1(\mathbb{R}^{d})}
$$
with an $x_0$-independent implicit constant.
\findem

\mb It remains us to prove the following lemma.

\begin{lem} \label{lemme} The bilinear map $(f,g) \mapsto f \otimes g$ is bounded from $\mathcal{H}^1(\mathbb{R}^d) \times \mathcal{H}^1(\mathbb{R}^d)$ to $\mathcal{H}^1(\mathbb{R}^{2d})$.
\end{lem}

\dem This fact can be easily proved using the atomic decomposition of Hardy spaces or using the characterization with maximal functions (see \cite{Stein}). Let us describe these two points of view.

\mb {\bf Use of maximal functions: } We recall that the Hardy spaces on $\R^d$ can be defined using the following norm~:
$$ \|f\|_{\mathcal{H}^1(\mathbb{R}^{d})} \simeq \left\| x\to \sup_{t} \left|\int_{\R^d} f(x-u) e^{-|u|^2/t^2} \frac{du}{t^{d}} \right| \right\|_{L^1(\R^d)}.$$
So we get easily~:
\begin{align*}
 \left\| f \otimes g \right\|_{\mathcal{H}^1(\mathbb{R}^{2d})} & \simeq \left\| (x,y)\to \sup_{t} \left|\int_{\R^{2d}} f(x-u)g(y-v) e^{-(|u|^2+|v|^2)/t^2} \frac{dudv}{t^{2d}} \right| \right\|_{L^1(\R^{2d})} \\
& \lesssim \left\| \left[x\to \sup_{t} \left|\int_{\R^d} f(x-u)e^{-|u|^2/t^2} \frac{du}{t^d}\right|\right]
\left[y \to \sup_{t>0} \left|\int_{\R^d} g(y-v) e^{-|v|^2/t^2} \frac{dv}{t^{d}} \right| \right] \right\|_{L^1(\R^{2d})} \\
& \lesssim \|f\|_{\mathcal{H}^1(\mathbb{R}^{d})} \|g\|_{\mathcal{H}^1(\mathbb{R}^{d})}.
\end{align*}

\mb {\bf Use of atomic decomposition: } It just suffices to prove that for $f$ and $g$ two atoms on $\R^d$ then $f \otimes g\in \mathcal{H}^1(\mathbb{R}^{2d})$. So let assume that $f$ is an atom corresponding to a ball $B_1$ and $g$ to a ball $B_2$ with $r_{B_2}\geq r_{B_1}$. We chose $(Q_i)_{i}$ a finite and bounded covering of $B_2$ with balls of radius equal to $r_{B_1}$ and then $(\phi_{B_i})_i$ a corresponding partition of unity. Then we write:
$$ f\otimes g = \sum_{i} f \otimes (\phi_{Q_i} g) \overset{def}{=} \sum_{i} \left(\frac{|Q_i|}{|B_2|}\right)^{1/2} \frac{\|\phi_{Q_i} g\|_{L^2}}{\|g\|_{L^2}} b_i.$$
As
$$ \sum_{i} \left(\frac{|Q_i|}{|B_2|}\right)^{1/2} \frac{\|\phi_{Q_i} g\|_{L^2}}{\|g\|_{L^2}} \leq \left(\sum_{i} \frac{|Q_i|}{|B_2|}\right)^{1/2} \left( \sum_{i} \frac{\|\phi_{Q_i} g\|_{L^2}^2}{\|g\|_{L^2}^2} \right)^{1/2} \lesssim 1,$$
it just suffices to check that $b_i$ is an atom for the ball $B_1\times Q_i$. The cancellation property for $f$ implies $\int b_i =0$ and we have
\begin{align*}
 \| b_i\|_{L^2(B_1\times Q_i)} & \leq \left(\frac{|B_2|}{|Q_i|}\right)^{1/2} \frac{\|g\|_{L^2}}{\|\phi_{Q_i} g\|_{L^2}} \|f\|_{L^2(B_1)} \|\phi_{Q_i}g\|_{L^2} \\
 & \lesssim \left(\frac{|B_2|}{|Q_i|}\right)^{1/2} |B_2|^{-1/2} |B_1|^{-1/2} \\
 & \lesssim |B_1\times Q_i|^{-1/2}.
\end{align*}
So $b_i$ is an atom to the ball $B_1 \times Q_i$ and then $f\otimes g$ belongs to $\mathcal{H}^1(\mathbb{R}^{2d})$.
\findem

\mb We emphasize that the previous estimate with Hardy spaces is weaker than the one involving the $L^1$ spaces. However Hardy spaces can be interpolated with Lebesgue spaces : we refer the reader to Subsection 4.2 of \cite{B} for bilinear real interpolation involving Lebesgue and Hardy spaces.

\section{Optimality of the estimates} \label{sec:5}

In this section, we want to prove that the set of Lebesgue exponents for which we proved estimates in theorems~\ref{thm:thm3smooth} and~\ref{thm:thm3cm} are optimal; we consider the case where $m$ is homogeneous of degree $0$ (or at least has 0-homogeneous bounds) and the phase $\phi$ is a homogeneous polynomial of degree 2.

\mb It follows then from the scaling that the only possible estimates are
\begin{equation}
\label{squirrel}
\left\|B_\lambda(f,g) \right\|_{L^r} \lesssim |\lambda|^{-\frac{d}{2}\left(\frac{1}{p}+\frac{1}{q}-\frac{1}{r}\right)} \|f\|_{L^p} \|g\|_{L^q}.
\end{equation}
We claim that the best possible set of exponents $p,q,r$ is given in Theorem~\ref{thm:thm3smooth} by (\ref{eq:condi}).

\begin{prop}
\label{pingouin}
Assume that $m$ is homogeneous of degree $0$ (or has corresponding bounds), $\phi$ is homogeneous of degree 2, and $\phi$ satisfies the nondegeneracy conditions of Theorem~\ref{thm:thm3smooth}. Then the estimate~\ref{squirrel} holds in general only for $p,q,r$ satisfying (\ref{eq:condi}).
\end{prop}

\dem To prove the proposition, it suffices of course to build up counterexamples. See the figure \ref{fig:tetra} for the admissible exponents; they form a tetrahedron with vertices.

We will prove in the next two subsections that if one sets $m$ identically equal to one, there exists a non-degenerate
quadratic polynomial $\phi$, and functions $f$ and $g$ such that
\begin{itemize}
\item (First counterexample) $f$ and $g$ belong to $L^p$ (so $p=q$), but $B_1(f,g)$ does not belong to $L^r$ if $r = \frac{2p}{2-p}$ and $p < \frac{6}{5}$. This permits to prove the necessity of the second equation in (\ref{eq:condi}): to see this graphically, the line $(p,q,r) = (p,p,\frac{2p}{2-p})$ is plotted below (dotted line) in the coordinates $(\frac{1}{p},\frac{1}{q},\frac{1}{r'})$.
\begin{figure}[htbp!]
 \centering
\includegraphics[width=0.25\textwidth]{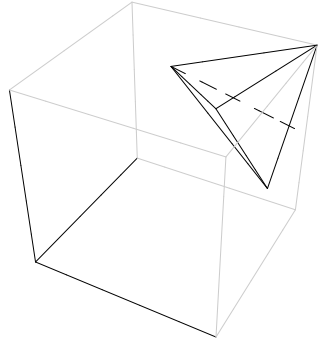}
\caption{The line $(p,q,r)$ = $(p,p,\frac{2p}{2-p})$ in the coordinates $\left(\frac{1}{p},\frac{1}{q},\frac{1}{r'} \right)$}
\label{fig:tetra4}
\end{figure}
The first counterexample stated above means that the statement of Theorem~\ref{thm:thm3smooth} becomes wrong as soon as $(p,q,r)$ are on the dotted line, but outside of the solid tetrahedron.

\item (Second counterexample) $f$ and $g$ belong to $L^p$ (so $p=q$), but $B_1(f,g)$ does not belong to $L^r$ if
$r=\frac{p}{p-1}$ and $p >  \frac{3}{2}$. This counterexample shows the necessity of the first condition in
(\ref{eq:condi}): to see this graphically, the line $(p,q,r) = (p,p,\frac{p}{p-1})$ is plotted below (dotted line) in the coordinates $(\frac{1}{p},\frac{1}{q},\frac{1}{r'})$.
\begin{figure}[htbp!]
 \centering
\includegraphics[width=0.25\textwidth]{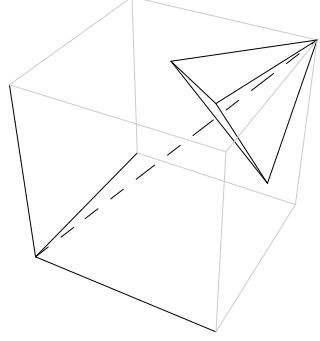}
\caption{The line $(p,q,r)$ = $(p,p,\frac{p}{p-1})$ in the coordinates $\left(\frac{1}{p},\frac{1}{q},\frac{1}{r'} \right)$}
\label{fig:tetra3}
\end{figure}
The second counterexample stated above means that the statement of Theorem~\ref{thm:thm3smooth} becomes wrong as soon as $(p,q,r)$ are on the dotted line, but outside of the solid tetrahedron.
\end{itemize}
By symmetry, we deduce similar results for the third and the fourth condition in (\ref{eq:condi}). This completes the proof of the proposition.
\findem

\subsection{The first counterexample}

Define the function $f$ on $\mathbb{R}^d$ by its Fourier transform
$$
\widehat{f}(\xi)\overset{def}{=} \frac{1}{\< |\xi| \> ^\mu}\;\;\;\;\;\mbox{where we denote $\<x\> = \sqrt{1 + x^2}$}.
$$
Take $m =1$, and $\phi(\xi,\eta) = \xi^2 + \eta^2 + 4\xi\cdot \eta$. The resulting operator $B_1$, when applied to $f$ and $f$, reads
$$
B(f,f)(x) = \int \int e^{ix\cdot(\eta+\xi) } e^{i(\xi^2 + \eta^2 + 8\xi\cdot \eta)} \frac{1}{\< |\xi| \> ^\mu} \frac{1}{\< |\eta| \> ^\mu} \,d\eta\, d\xi.
$$
We will prove the

\begin{lem}
As $|x|\rightarrow \infty$,
$$
|B(f,f)(x)| \sim \frac{C}{|x|^{2\mu}}.
$$
\end{lem}

\mb A simple computation shows that this provides the first counterexample needed in the proof of Proposition~\ref{pingouin}: indeed, it is easy to see that $f$ decays fast at infinity, and is smooth except at zero, where it has a singularity like $\frac{1}{|x|^{d-\mu}}$. Thus, if $\mu=d-\frac{d}{p}+\epsilon$, with $\epsilon>0$, $f$ belongs to $L^p$, and the lemma implies that $B(f,f)$ does not belong to $L^r$, for $r < \frac{d}{2\mu} = \frac{1}{2-\frac{2}{p}+2\frac{\epsilon}{d}}$. In particular, if $p<\frac{6}{5}$, $B(f,f)$ does not belong to $L^{\frac{2p}{2-p}}$ for $\epsilon$ small enough. \\
So it remains to prove the Lemma.

\dem It is convenient to change variables
$$
\xi'=\frac{\xi}{|x|} \;\;\;,\;\;\;\eta' = \frac{\eta}{|x|} \;\;\;,\;\;\;\omega = \frac{x}{|x|}\;\;\;\mbox{and}\;\;\;\psi(\omega,\xi',\eta') = \omega\cdot (\xi'+\eta') + \xi'^2 + \eta'^2 + 8\xi' \cdot \eta'
$$
and write
$$
B(f,f)(x) = |x|^{2d-2\mu} \int \int e^{i|x|^2 \psi(\omega,\eta',\xi') }  \frac{1}{|x|^{-\mu} \< |\xi'||x| \> ^\mu}  \frac{1}{|x|^{-\mu} \< |\eta'||x| \> ^\mu} \,d\eta'\, d\xi'.
$$
As usual, this integral will be estimated by splitting it into several pieces corresponding to different integration domains, and estimating them separately. Since we already went through this procedure a number of times, we will be a bit sketchy, and in particular not write the cut-off functions.

\bigskip
\noindent
\underline{Estimate of the piece $|(\xi',\eta')|\sim 1$} This corresponds to
$$
G_1 (x) \overset{def}{=} |x|^{2d-2\mu} \int \int_{|(\xi',\eta')|\sim 1} e^{i|x|^2 \psi(\omega,\eta',\xi') }
\frac{1}{|x|^{-\mu} \< |\xi'||x| \> ^\mu}  \frac{1}{|x|^{-\mu} \< |\eta'||x| \> ^\mu}\,d\eta'\, d\xi'.
$$
This is the domain where the phase is stationary: indeed, the derivative in $\eta',\xi'$ of $\psi(\omega,\xi',\eta')$ vanishes for $\xi'=\eta'=-\frac{\omega}{10}$. Thus the stationary phase lemma gives
$$
\mbox{as $x\rightarrow \infty$}, \;\;\;\;\;|G_1(x)| \sim \frac{C}{|x|^{2\mu}}.
$$

\bigskip
\noindent
\underline{Estimate of the piece $|(\xi',\eta')| >> 1$} This is
$$
G_2 (x) \overset{def}{=} |x|^{2d-2\mu} \int \int_{|(\xi',\eta')|>> 1} e^{i|x|^2 \psi(\omega,\eta',\xi') }  \frac{1}{|x|^{-\mu} \< |\xi'||x| \> ^\mu}  \frac{1}{|x|^{-\mu} \< |\eta'||x| \> ^\mu} \,d\eta'\, d\xi'.
$$
On this integration domain, the phase satisfies $|\nabla_{\eta',\xi'} \psi| \gtrsim |(\eta',\xi')|$. Repeated integration by parts give
$$
|G_2(x)| \lesssim |x|^{-N} \;\;\;\;\;\mbox{for any $N$}.
$$

\bigskip
\noindent
\underline{Estimate of the piece $ |(\xi',\eta')| << 1$} We are now considering
$$
G_3(x) \overset{def}{=} |x|^{2d-2\mu} \int \int_{\frac{1}{|x|} \lesssim |\xi',\eta'| << 1  } e^{i|x|^2 \psi(\eta',\xi') }  \frac{1}{|x|^{-\mu} \< |\xi'||x| \> ^\mu}  \frac{1}{|x|^{-\mu} \< |\eta'||x| \> ^\mu} \,d\eta'\, d\xi'.
$$
On this integration domain, the phase satisfies $|\nabla_{\eta',\xi'} \psi| \gtrsim 1$. Once again, repeated integration by parts give
$$
|G_3(x)| \lesssim |x|^{-N} \;\;\;\;\;\mbox{for any $N$}.
$$

\gb Since $B(f,f)=G_1+G_2+G_3$, the lemma follows from the above estimates. \findem

\subsection{The second counterexample} Set first
$$
F(x) = \frac{\phi(x)}{|x|^{d-\mu}},
$$
with $\phi$ a smooth function such that $\phi = 1$ on $B(0,\frac{1}{2})$, and $\phi=0$ on $B(0,1)^c$. It is easy to see with the help of the stationary phase lemma that $e^{it\Delta} F$ is a smooth function such that, for fixed $t$,
$$
\mbox{as $x \rightarrow \infty$}, \;\;\;\;e^{it\Delta} F \sim C_t \frac{e^{i\frac{t}{4}x^2}}{|x|^\mu}.
$$
Then set
$$
f \overset{def}{=} e^{-8i\Delta} F \;\;\;\;\mbox{and}\;\;\;\; g \overset{def}{=} e^{-i \Delta} F
$$
and
$$
\phi(\eta,\xi) \overset{def}{=} 8 \xi^2 + \eta^2 + (\xi+\eta)^2\;\;\;\;\mbox{and}\;\;\;\;m\overset{def}{=}1.
$$
Then
$$
B_\lambda(f,g) = e^{i\Delta}\left( e^{8i\Delta} f e^{i\Delta} g \right) = e^{i\Delta} F^2 = e^{i\Delta} \frac{\phi(x)^2}{|x|^{2d-2\mu}} \sim \frac{e^{i\frac{t}{4}x^2}}{|x|^{2\mu-d}}
$$
as $x$ goes to infinity. Picking $\mu = \frac{d}{p}+\epsilon$, one sees that $f$ and $g$ belong to $L^p$, but that $B_\lambda(f,g)$ does not belong to $L^{\frac{p}{p-1}}$ for $p>\frac{3}{2}$.

\section{Applications - Boundedness of new bilinear multipliers with non smooth symbols} \label{sec:6}

In this section we want to describe an application of these estimates. Mainly we give boundedness results on Lebesgue spaces for new singular bilinear multipliers (belonging to no known classes). These new singular symbols will be defined using the notion of ``finite part'', which we shall first make precise.

We mention that another class of singular symbols was considered by Kenig and Stein~\cite{KS}.

\subsection{``Finite part'' of the inverse of a smooth function}

Let us temporarily forget our subject of bilinear oscillatory integrals and consider a smooth function $\phi$ on $\R^n$. It is obvious that the inverse function $1/\phi$ could be locally non-integrable around the domain $\{\phi=0\}$ and so it is not clear when $1/\phi$ could define a distribution. \\
In order to get around this problem, many works have dealt in some particular cases with the ``principal value'' of $1/x$, or the finite part of $1/x^2$ on $\R$.

\mb Here, we would like to use oscillatory integral in order to give a precise sense to a distribution, which we call ``Finite part of $1/\phi$'' and then to give some examples in Subsection \ref{subsec:ex}.\\
We assume the following~: for some exponent $\delta>1$ and any real $\lambda\geq 1$, we have \be{decay} \left|
\int_{\R^{n}} e^{i\lambda\phi(x)} f(x) dx \right| \lesssim \lambda^{-\delta} C(f), \ee where $C(f)$ is one of the
semi-norms of the Schwartz space $\s(\R^n)$.

\begin{thm} \label{thm:def}
Let $\phi$ a real smooth function such that its derivatives admit an at most polynomial growth and satisfying (\ref{decay}). Then the sequence of distributions
 $$ \left(-\frac{e^{iT\phi}-1}{\phi}\right)_{T\geq 1}$$
has a limit in $\s'(\R^n)$ for $T\to\infty$. We call it ``Finite part of $1/\phi$'' and it will be denoted by $F.P. (\frac{1}{\phi})$.
\end{thm}

\dem First, for each $T\geq 1$, $D_T\overset{def}{=}-\frac{e^{iT\phi}-1}{\phi}$ is bounded in $L^\infty$ by the constant $T$ so it is really a distribution. For any test function $f\in\s(\R^n)$, we have~:
$$ \langle D_T,f \rangle \overset{def}{=} -\int_{\R^n} \frac{e^{iT\phi(x)}-1}{\phi(x)} f(x) dx =  -\int_{\R^n} \frac{e^{i(T+1)\phi(x)}-e^{i\phi(x)}}{\phi(x)} g(x) dx$$
with $g(x)\overset{def}{=}e^{-i\phi(x)} f(x)$. Since $\phi$ has at most a polynomial behavior, $g$ still belongs to $\s(\R^d)$. We use the following formula~:
$$ \langle D_T,f \rangle = -i\int_{1}^T \int_{\R^n} e^{i\lambda \phi(x)} g(x) dx d\lambda, $$
which gives us (according to (\ref{decay})) for $T<T'$:
$$ \left| \langle D_T-D_{T'},f \rangle \right| \lesssim \int_{T}^{T'} \lambda^{-\delta} C(g) d\lambda \lesssim T^{-\delta+1}C'(f)$$
with an other norm $C'$. As $\delta>1$, we deduce that $(D_T)_T$ is a Cauchy sequence in the space $\s'(\R^n)$ and so converges. \findem

\mb By this way, we have defined the finite part distribution as an abstract limit. We refer the reader to the next subsection for examples.

\mb This procedure can be seen as follows : when the function $\phi$ too fastly decreases around $\{\phi=0\}$ then the function $1/\phi$ is not locally integrable. So we have to regularize this operation, to define a distribution. Aiming that, we define the finite part distribution, which corresponds to remove the mean of the function on the characteristic domain $\{\phi=0\}$ and to create some extra-cancellation around this singular region. A technical difficulty is to give a precise sense to the mean function on such a set (which is not assumed to be smooth).\\
We study in the following subsection two cases, where we explain this point of view (Propositions \ref{prop:1}, \ref{prop:2} and \ref{prop:3}) and then we come back to our bilinear oscillatory integrals in Subsection \ref{subsec:bilinear}.

\subsection{Particular cases of ``Finite part'' distributions} \label{subsec:ex}

We would like to compute the limit of the distribution
$$D_T\overset{def}{=}  -\frac{e^{iT\phi}-1}{\phi}$$ as $T\to\infty$.
We have
$$ \langle D_T,f\rangle\overset{def}{=} - \int_{\R^n} \left[\frac{e^{iT\phi(x)}-1}{\phi(x)} \right] f(x) dx.$$
First for any $\epsilon>0$ and smooth function $f$, we have~:
$$ \lim_{T\to \infty}  \ \int_{|\phi|\geq \epsilon}  \left[\frac{e^{iT\phi(x)}-1}{\phi(x)} \right] f(x) dx = - \int_{|\phi|\geq \epsilon} \frac{1}{\phi(x)} f(x) dx.
$$
So it remains to study the behavior around the manifold $\Delta\overset{def}{=}\{\phi=0\}$ and we deal with
$$\langle D^1_T,f\rangle \overset{def}{=}  -\int_{|\phi|\leq \epsilon} \left[ \frac{e^{iT\phi(x)}-1}{\phi(x)} \right] f(x) dx$$
for a small enough $\epsilon>0$.

\gb
{\bf First case :} no critical points on $\Delta\overset{def}{=}\phi^{-1}(0)$.

\mb We assume that $\phi$ has no critical points so $|\nabla \phi|$ is non vanishing on $\phi^{-1}(0)$. In this case, we know that $\Delta=\phi^{-1}(0)$ is a smooth hypersurface of $\R^{2d}$.
Then the ``level-set function integration formula'' gives~: for any continuous function $h$ whose support is close enough to $\Delta$,
$$ \int h(x) dx = \int_{-1}^1 \int_{\{\phi=t\}} h(x) \frac{d\sigma_t(x)}{|\nabla\phi(x)|} dt$$
where $d\sigma_t$ is the superficial measure on $\phi^{-1}(t)$ (which is a smooth hypersurface for small enough $t$).
So let us write $D^1_T$ as follows~:
$$\langle D^1_T,f \rangle \overset{def}{=} -\int_{-\epsilon}^\epsilon \left[\frac{e^{iTt}-1}{t} \right] \zeta(t) dt$$
with
$$ \zeta(t)\overset{def}{=} \int_{\{\phi=t\}} f(x) \frac{d\sigma_t(x)}{|\nabla\phi(x)|}.$$
Using that
\be{eq:h} \lim_{T\to \infty} p.v. \int_{-\epsilon}^\epsilon \frac{e^{iTt}}{t} \zeta(t) dt =\widehat{p.v. \frac{1}{\cdot}} (1) \zeta(0)=i\pi \zeta(0), \ee
 we finally deduce
$$ \lim_{T\to\infty} \langle D^1_T,f\rangle = -i\pi\zeta(0)+p.v.\int_{-\epsilon}^\epsilon \frac{1}{t} \zeta(t) dt$$
So using the same manipulation with the level-sets, we finally get
$$ \langle F.P. (1/\phi),f\rangle = \lim_{T\to \infty} \langle D_T,f\rangle =  -i\pi\zeta(0) + p.v. \int \frac{1}{\phi(x)} f(x) dx,$$
where the principal value is associated to the integral in the $t$-variable.

\begin{prop} \label{prop:1} In this case, we have the following result~:
$$  F.P.(1/\phi) =  p.v. \frac{1}{\phi} - i\pi \frac{d\sigma_\Delta}{|\nabla\phi|},$$
where the principal value is in the sense of the level set for $\phi$~:
$$ p.v. \frac{1}{\phi} = \lim_{r\to 0} \frac{1}{\phi} {\bf 1}_{|\phi|>r}$$
and $\sigma_\Delta$ is the superficial measure of $\Delta=\phi^{-1}(0)$.
\end{prop}

\gb
{\bf Second case :} Locally around a critical point.

\mb We assume now that we work locally around $x_0$ a non-degenerate critical point of $\phi$ belonging to $\phi^{-1}(0)$. Using the assumptions of non-degenerescence, we know that these points are separated.

\mb Then around each critical point, we can repeat the previous arguments, we do not know if $\Delta$ is a smooth hypersurface or not, however the ``level-set function integration formula'' still holds as the gradient is non vanishing around $x_0$. Indeed the set $\phi^{-1}(0)$ has a vanishing $n$-dimensional measure as Morse's Lemma reduces $\phi^{-1}(0)$ to some quadric-surface in an appropriate system of coordinates. So we can exactly apply the same reasonning in order to obtain~:
$$ \langle D^1_T,f \rangle \overset{def}{=} -p.v. \int_{-\epsilon}^\epsilon \left[\frac{e^{iTt}-1}{t} \right] \zeta(t) dt$$
with
$$ \zeta(t)\overset{def}{=} \int_{\{\phi=t\}} f(x) \frac{d\sigma_t(x)}{|\nabla\phi(x)|}.$$
The new problem here is that $\zeta(0)$ is non-defined as $\nabla\phi=0$ on $x_0\in\phi^{-1}(0)$. Thus one has to be careful. Due to the non-degenerescence, we know that locally around $x_0$~:
$$ |\nabla \phi(x)| \simeq |x-x_0|.$$

\mb So there is an extra difficulty in order to estimate the behavior of $\zeta(t)$ when $t\to 0$. Let us explain some particular situations for $n\geq 3$, depending on the signature of the Hessian matrix.

\mb $\bigstar$ For example, assume that the Hessian $\operatorname{Hess}\phi(x_0)$ is strictly positive. Then we know
that around $x_0$
 $$\phi(x) \simeq |x-x_0|^2.$$
 So it is obvious that for $n> 2$ the function~: $$ \frac{1}{|\phi(x)|} \lesssim \frac{1}{|x-x_0|^2}$$
 is integrable around the critical point $x_0\in\R^{n}$. So it is not necessary to define a principal value of $1/\phi$ and we deduce the next result.

\begin{prop} \label{prop:2} For $n\geq 3$, if $\phi$ has one non-degenerate critical point $x_0$ with a positive (or negative) hessian matrix, then
$$ F.P. (1/\phi) = \frac{1}{\phi} \in L^1_{loc}.$$
\end{prop}

\mb $\bigstar$ We claim that in the other cases the function $1/\phi$ is not integrable around $x_0$ and so we need to invoke a ``finite part'' in order to give a distributional sense to $1/\phi$. So consider that the hessian matrix has a signature $(p,n-p)$ with $1\leq p \leq n-1$. Then as previously, according to Morse's Lemma we can compare
 $$ \frac{1}{|\phi(x)|} \lesssim \frac{1}{\left|\left(|y_1|^2+..|y_p|^2\right)-\left(|y_{p+1}|^2+..+|y_{n}|^2 \right)\right|},$$
where $y\in\R^{n}$ is a new system of coordinates (around $0$ when $x$ is around $x_0$) with $y_1,..,y_{p}$ corresponding to the positive eigenvalues of $\operatorname{Hess}\phi(x_0)$ and $y_{p+1},..,y_{n}$ corresponding to the negative ones. So we have to consider the quadric surface
$$ S\overset{def}{=} \left\{y, |y_1|^2+..|y_p|^2 = |y_{p+1}|^2+..+|y_{n}|^2 \right\}.$$
 We note that $S$ is a hypersurface of dimension $n-1$ due to $p\neq 0$ and $p\neq n$ (if $p\in\{0,n\}$, $S$ is reduced to one point in $\R^{n}$). Then obviously, we have with $U$ any small enough neighborhood around $x_0$
 $$ \int_U \frac{1}{|\phi(x)|} dx \gtrsim \int_0^\epsilon \int_0^r \frac{1}{r^2-t^2} t^{p-1} r^{n-p-1} dt dr = \infty.$$
 So when the hessian matrix $\operatorname{Hess}\phi(x_0)$ is non-positive and non-negative, $1/\phi$ is non integrable locally around $x_0$, and this is why one has to define a ``finite part''.
Up to a change of variables, Morse's Lemma implies that we can assume
$$ \phi(x)=\left(|x_1|^2+..+|x_p|^2\right)-\left(|x_{p+1}|^2+..+|x_{n}|^2 \right).$$
So we get~:
$$ \langle D_T,f\rangle = -\int \frac{e^{iT\phi(x)} -1}{\phi(x)} f(x) dx = -\int_{r} \int_{t} \frac{e^{iT(t^2-r^2)}}{t^2-r^2} t^{p-1}r^{n-p-1}  \zeta(t,r)dtdr,$$
with
$$\zeta(t,r) \overset{def}{=} \int_{\left\{\genfrac{}{}{0pt}{}{|x_1|^2+..|x_p|^2=t^2}{|x_{p+1}|^2+..+|x_{n}|^2=r^2}\right\}} f(x) dx.$$
We leave the details to the reader, as previously using (\ref{eq:h}), we obtain~:
$$ \lim_{T\to \infty} \int \frac{e^{iT\phi(x)}}{\phi(x)} f(x) dx = i\pi \int r^{n-3} \zeta(r,r) dr.$$
Consequently, we get the following result~:
$$F.P. (1/\phi) = -i\pi d\sigma +  p.v. \frac{1}{\phi},$$
with $d\sigma$ defined as
$$ d\sigma= r^{n-3} d\sigma_{\left\{\genfrac{}{}{0pt}{}{|x_1|^2+..|x_p|^2=r^2}{|x_{p+1}|^2+..+|x_{n}|^2=r^2}\right\}}(x) dr = d\sigma_S.$$
Using the change of variables given by Morse's Lemma, we finally get the following result~:

\begin{prop} \label{prop:3} For $n\geq 3$, if $\phi$ has one non-degenerate critical point $x_0$ with a non positive and non negative Hessian matrix, then we have
$$  F.P.(1/\phi) =  p.v.\frac{1}{\phi} - i\pi \frac{d\sigma_\Delta}{|\nabla\phi|},$$
with the same notations as for Proposition \ref{prop:1}.
\end{prop}

\mb We remark that in this particular case, we know that $|\nabla\phi(x)| \simeq |x-x_0|$ so $|\nabla\phi|^{-1}$ is integrable along the manifold $\phi^{-1}(0)$ as its dimension is equal to $n-1\geq 2>1$.

\mb We note that both Proposition \ref{prop:3} and Proposition \ref{prop:1} express the same result.

\begin{rem} The distribution $\frac{d\sigma_\Delta}{|\nabla\phi|}$ can be seen as the distribution image of the dirac distribution $\delta_0$ via the map $\phi$, so in a {\it certain sense}, we can write
 $$\frac{d\sigma_\Delta}{|\nabla\phi|} =\delta_0(\phi).$$
\end{rem}

\begin{rem} It is interesting to note that the non-degeneracy assumption on the hessian matrix corresponds to a non-vanishing curvature on the manifold $\phi^{-1}(0)$. Indeed we know that the principal curvatures are fixed by the eigenvalues of the restriction of the Hessian matrix to the tangent space, and so are bounded below.
\end{rem}

\subsection{Boundedness of singular bilinear multipliers} \label{subsec:bilinear}

Let us now come back to our main purpose. So consider a smooth function $\phi(\eta,\xi)$ on $\R^{2d}$. We assume that we can define in the previous sense a finite part of $1/\phi$ in $\s'(\R^{2d})$. Suppose moreover that we have a decay for some exponents $p,q,r\in[1,\infty)$ and any $|\lambda|\geq 1$
\be{eq:decay} \left\|B_\lambda(f,g)\right\|_{L^r} \lesssim |\lambda|^{-\rho} \|f\|_{L^p}\|g\|_{L^q} \ee
for some exponent $\rho>1$. \\
We recall our bilinear oscillatory integral~:
$$B_\lambda(f,g)(x)\overset{def}{=} \frac{1}{(2\pi)^{d/2}}\int_{\R^{2d}} e^{ix\cdot(\eta+\xi)} e^{i\lambda \phi(\eta,\xi)} m(\eta,\xi) \widehat{f}(\eta)\widehat{g}(\xi) d\eta d\xi.$$

\begin{rem} We note that the power $\rho$ obtained in Theorem \ref{thm:thm3smooth} is smaller (or equal) than $1$ if $d=1$. So we have to consider multidimensional variables ($d\geq 2$).
\end{rem}

\mb Then as done for Theorem \ref{thm:def}, we have also the following bilinear version~:

\begin{thm} \label{thm:newmul} Under the above assumptions, the bilinear multiplier associated to the symbol
$$
e^{i\phi(\eta,\xi)} F.P. \frac{1}{\phi(\eta,\xi)} m(\eta,\xi)$$
is bounded from $L^p(\R^d) \times L^q(\R^d)$ into $L^r(\R^d)$.
\end{thm}

\mb According to the examples of ``Finite Parts'' given by Propositions \ref{prop:1} and \ref{prop:3}, it appears singular symbols supported on some submanifold on $\R^{2d}$. We move the reader to Section \ref{sec:singular} for some bilinear oscillatory integrals involving such symbols.

\mb We want to finish this section by giving an example.

\begin{ex}
Let us consider the phase $\phi(\eta,\xi)=\xi.\eta$ and $m=1$. Then Theorem \ref{thm:thm3cm2} with $\lambda=1$ implies that the bilinear operator $T_{\sigma}$ with $\sigma=e^{i\phi}$ 
belongs to the class ${\mathcal M}_{p,q,r}$ for exponents $p,q,r$ in the appropriate range.

Then Theorem \ref{thm:newmul} implies that the bilinear multiplier associated to the symbol
$$  F.P. \frac{e^{i\eta \cdot \xi}}{\eta.\xi} $$
belongs to ${\mathcal M}_{p,q,r}$ as soon as the additional condition
$$ \frac{d}{2}\left(\frac{1}{p}+\frac{1}{q}-\frac{1}{r}\right)>1.$$
is satisfied.
\end{ex}

\subsection{Application to scattering theory}

Consider the PDE introduced at the beginning of this article
$$ 
\left\{ \begin{array}{l}
\partial_t u + iP(D) u = T_m(u,u) \\
u(t=0)=u_0.
\end{array} \right.
$$
The solution $u$ is said to scatter if
$$ f(t,x) = e^{itP(D)}[u(t,.)](x).$$
has a limit $f_\infty$ (for some topology) as $t$ goes to infinity. 

It is possible to solve this equation iteratively by setting $u^1 = e^{itP(D)}u_0$ and, for any $n \geq 1$, defining $u^n$ by
$$
\left\{ \begin{array}{l}
\partial_t u^n + iP(D) u^n = T_m(u^{n-1},u^{n-1}) \\
u^n(t=0)=0
\end{array} \right.
$$
(thus, $u^n$ is an $n$-linear operator in $u_0$).
Ignoring all convergence questions, we get $u = \sum_{n=1}^\infty u^n$, which means, still working formally,
$$
f_\infty = \sum_{n=1}^\infty L^n(u_0),\;\;\mbox{where $L^n(u_0) = \lim_{t \rightarrow \infty} e^{itP(D)}[u^n(t,.)]$}
$$
is an $n$-linear operator. An easy computation using Duhamel's formula gives
$$
L^1(u_0) = u^0 \;\;\;\;\mbox{and}\;\;\;\;L^2(u_0)= B_\mu(u_0,u_0)\;\;\;\;\mbox{with}\;\;\;\;\mu(\xi,\eta) = \frac{e^{i[P(\xi+\eta) - P(\eta) - P(\xi)]}}{P(\xi+\eta) - P(\eta) - P(\xi)};
$$
in other words, the second derivative of the scattering operator is given by a principal value operator similar to the ones which have just been studied.

\end{document}